\documentclass[11pt]{amsart}
\usepackage{amsmath,amsthm,amssymb,url}
\usepackage[all]{xy}

\newtheorem*{theorem*}{Theorem}
\newtheorem*{proposition*}{Proposition}

\newtheorem{theorem}{Theorem}[section]
\newtheorem{lemma}[theorem]{Lemma}
\newtheorem{proposition}[theorem]{Proposition}

\theoremstyle{definition}
\newtheorem{definition}[theorem]{Definition}

\theoremstyle{remark}
\newtheorem{remark}[theorem]{Remark}

\newcommand{\mapor}[1]{{\stackrel{#1}{\longrightarrow}}}

\renewcommand{\bar}{\overline}

\newcommand{\g}{\mathfrak{g}}
\newcommand{\m}{\mathfrak{m}}

\newcommand{\K}{\mathbb{K}}

\newcommand{\Id}{\operatorname{Id}}

\newcommand{\MC}{\operatorname{MC}}
\newcommand{\Def}{\operatorname{Def}}
\newcommand{\Del}{\operatorname{Del}}

\newcommand{\Spec}{\operatorname{Spec}}
\newcommand{\Ext}{\operatorname{Ext}}
\newcommand{\holim}{\mathop{\rm holim}}
\newcommand{\mathscr}{\mathbb}

\newcommand{\Art}{\mathbf{Art}}
\newcommand{\Set}{\mathbf{Set}}

\newcommand{\Tot}{\operatorname{Tot}}
\newcommand{\N}{\mathbb{N}}

\newcommand{\deltabar}{\bar{\partial}}
\newcommand{\Eps}{\mathcal{E}}
\newcommand{\U}{\mathcal{U}}

\renewcommand{\O}{\mathcal{O}}

\newenvironment{acknowledgement}{\par\addvspace{17pt}\small\rm
\trivlist\item[\hskip\labelsep{\it Acknowledgement.}]}
{\endtrivlist\addvspace{6pt}}

\begin{document}
\title[Infinitesimal deformations of coherent sheaves
]{Differential graded Lie algebras controlling infinitesimal
deformations of coherent sheaves}

\begin{abstract}
We use the Thom-Whitney construction to show that infinitesimal
deformations of a coherent sheaf ${\mathcal F}$ are controlled by
the differential graded Lie algebra of global sections of an
acyclic resolution of the sheaf $\Eps nd^*(\Eps^\cdot)$, where
$\Eps^\cdot$ is any locally free resolution of ${\mathcal F}$. In
particular, one recovers the well known fact that the tangent
space to $\Def_{\mathcal F}$ is $\Ext^1({\mathcal F},{\mathcal F})$, and
obstructions are contained in $\Ext^2({\mathcal F},{\mathcal F})$.
\par
The main tool is the identification of the deformation functor
associated with the Thom-Whitney DGLA of a semicosimplicial DGLA ${\mathfrak g}^\Delta$,
whose cohomology is concentrated in nonnegative degrees, with a
noncommutative \v{C}ech cohomology-type functor $H^1_{\rm sc}(\exp
{\mathfrak g}^\Delta)$.
\end{abstract}


\date{}
\author{Domenico Fiorenza}
\address{\newline Dipartimento di Matematica \lq\lq Guido
Castelnuovo\rq\rq,\hfill\newline Sapienza Universit\`a di Roma,
\hfill\newline P.le Aldo Moro 5, I-00185 Roma Italy.}
\email{fiorenza@mat.uniroma1.it}
\urladdr{www.mat.uniroma1.it/\~{}fiorenza/}

\author{Donatella Iacono}
\address{\newline Max-Planck Institut f\"ur Mathematik,\hfill\newline Vivatsgasse 7, 
D 53111 Bonn Germany}
\email{iacono@mpim-bonn.mpg.de}

\author{Elena Martinengo}
\address{\newline Dipartimento di Matematica \lq\lq Guido
Castelnuovo\rq\rq,\hfill\newline Sapienza Universit\`a di Roma,
\hfill\newline P.le Aldo Moro 5, I-00185 Roma Italy.}
\email{martinengo@mat.uniroma1.it}
\urladdr{www.mat.uniroma1.it/dottorato/}

\newpage
\maketitle
\section*{Introduction}

The classical approach to deformation theory, starting  with
Kodaira and Spencer's studies on deformations of complex
manifolds, consists in deforming the objects locally and then glue
back together these local deformations. During the last thirty
years, another approach to deformation problems has been
developed. The philosophy underlying it, essentially due to
Quillen, Deligne, Drinfeld and Kontsevich, is that, in
characteristic zero, every deformation problem is controlled  by a
differential graded Lie algebra, via solutions of Maurer-Cartan
equation modulo gauge equivalence. The aim of this paper is to
exhibit an explicit equivalence between the two approaches for the
problem of infinitesimal deformations of coherent sheaves.

\bigskip

In the particular case of a locally free sheaf $\Eps$ of
${\mathcal O}_X$-modules on a complex manifold $X$, the
Kodaira-Spencer's description of deformations of $\Eps$ is given
in terms of the \v{C}ech functor $H^1(X; \exp \Eps nd(\Eps))$,
where $\Eps nd(\Eps)$ is the sheaf of endomorphism of $\Eps$.
Indeed, a  locally free sheaf has only trivial local deformations
and so a deformation of $\Eps$ is reduced to a deformation of the
gluing data of its local charts, and the compatibility conditions
these gluing data have to satisfy is precisely expressed by the
cocycle condition in the \v{C}ech functor. On the other hand, it
is well known that deformations of $\Eps$ are controlled by the
DGLA of global sections of an acyclic resolution of $\Eps nd(
\Eps)$, e.g., by the DGLA $A^{0,*}_X(\Eps nd (\Eps))$ of
$(0,*)$-forms on $X$ with values in the sheaf of endomorphisms of
the sheaf $\Eps$.

\bigskip
The equivalence between these two descriptions  is best understood
by moving from set-valued to groupoid-valued deformation functors;
see, e.g., \cite{hinich,Pridham}. Associating with any open set
$U$ in $X$ the groupoid $\Def_{\Eps|_U}$ of infinitesimal
deformations of $\Eps$ over $U$ (over a fixed base $\Spec A$, for
some local Artin ring $A$) defines a stack over ${\bf{Top}}_X$; this is
just a one-word way of saying that global deformations of $\Eps$
are the same thing as the descent data for its local deformations:
\[ \Def_\Eps \simeq \displaystyle \holim_{U\in \Delta_\U} \Def_{\Eps|_U}, \]
where $\Delta_\U$ is the semisimplicial object in ${\bf{Top}}_X$
associated with an open cover $\U$ of $X$. Next, one sees that
locally the groupoid of deformations of $\Eps|_U$ is equivalent to
the Deligne groupoid of $\Eps nd(\Eps)(U)$; since these
equivalences are compatible with restriction maps, one has an
equivalence of semicosimplicial groupoids. Finally, Deligne
groupoid commutes with homotopy limits of DGLA concentrated in
positive degree  (see \cite{hinich}), so that
\[  \Def_\Eps \simeq  \holim_{U\in \Delta_\U} \Del_{\Eps nd(\Eps)(U)}\simeq\Del_{\substack{{\holim \Eps nd(\Eps)(U)}\\{\scriptscriptstyle{U\in \Delta_\U\phantom{mmmmmi}}}}}.\]
This shows that the problem of infinitesimal  deformations of
$\Eps$ is controlled by the DGLA $\holim_{U\in\Delta_\U} \Eps nd(\Eps)(U)$. It is now a simple exercise in homological algebra showing
that there is a quasi-isomorphism of DGLAs
\[
\holim_{U\in\Delta_\U} \Eps nd(\Eps)(U)\simeq A^{0,*}_X(\Eps nd(\Eps)).
\]
The reader who prefers  to  not leave  the peaceful realm of
set-valued deformation functors can found a direct (but less
enlightening) proof of the equivalence between the
Kodaira-Spencer's and the DGLA approach to infinitesimal
deformation of locally free sheaves in  \cite{Fio-Man-Mart}, where
the explicit Thom-Whitney model for $\holim_{U\in\Delta_\U} \Eps
nd(\Eps)(U)$ is used.

\bigskip

We now turn our attention to deformations  of a coherent sheaf
$\mathcal F$ of $\O_X$-modules on a complex manifold or an
algebraic variety $X$. The classical approach to this deformation
problem is based on a locally free resolution $\Eps^\cdot\to
\mathcal F$ of $\mathcal F$; then, the data of a deformation of
$\mathcal F$ are the data of local deformations of $\Eps^\cdot$
with appropriate gluing conditions. More precisely, the sheaf of
differential graded Lie algebras $\Eps nd^*(\Eps^\cdot)$ of the
endomorphisms of the resolution $\Eps^\cdot$ controls
infinitesimal deformations of $\mathcal F$ via the \v{C}ech-type
functor $H^1_{\rm Ho}(X;\exp \Eps nd^*(\Eps^\cdot))$; the
subscript ${\rm Ho}$ refers to the fact that cocycle conditions
hold only up to homotopy. The functor $H^1_{\rm Ho}(X;\exp \Eps
nd^*(\Eps^\cdot))$ is actually independent of the particular
resolution chosen. And again, on the DGLA side, one proves that
infinitesimal deformations of ${\mathcal F}$ are controlled by the
DGLA of global sections of an acyclic resolution of $\Eps
nd^*(\Eps^\cdot)$; in particular, one recovers the well known fact
that the tangent space to $\Def_{\mathcal F}$ is
$\Ext^1({\mathcal F},{\mathcal F})$, and obstructions are contained in
$\Ext^2({\mathcal F},{\mathcal F})$.

\bigskip

To see why such a result should hold, one has to  make a further
step and go from groupoid-valued to $\infty$-groupoid-valued
deformation functors, and to think the whole problem in terms of
$\infty$-stacks \cite{hirschowitz-simpson, lurie, toen}. Indeed,
due to the presence of negative degree components in $\Eps nd^*(
\Eps^\cdot)$, the groupoids $\Def_{\mathcal F|_U}$ are no more
equivalent to the Deligne groupoids $\Del_{\Eps nd^*(\Eps^\cdot)(U)}$;
yet from the $\infty$-groupoid point of view it is natural to
expect that the stack $\Def_{\mathcal F}$ is locally homotopy
equivalent to the $\infty$-stack $\MC_{\bullet}(\Eps nd^*(
\Eps^\cdot))$. Then one reasons as in the locally free sheaf case,
using the fact that the Kan complexes-valued functor $\MC_\bullet$
commutes with homotopy limits of DGLAs whose cohomology is
concentrated in positive degree \cite{getzler}:
\[  \Def_{\mathcal F} \simeq  \holim_{U\in \Delta_\U} \Def_{\mathcal F|_U}\simeq  \holim_{U\in \Delta_\U} {\MC_{\bullet}(\Eps nd^*(\Eps^\cdot)(U))}\simeq {\MC_\bullet}(\holim_{U\in\Delta_\U} \Eps nd^*(\Eps^\cdot)(U)).\]
As above,  the homotopy limit $\holim_{U\in\Delta_\U}\Eps nd^*(\Eps^\cdot) (U)$ is quasiisomorphic to the DGLA of global sections
of an acyclic resolution of $\Eps nd^*(\Eps^\cdot)$, which therefore
controls the infinitesimal deformations of ${\mathcal F}$.

\bigskip

The aim of this paper is to give a direct proof  of this fact at
the level of set-valued deformation functors. The proof closely
follows the argument in \cite{Fio-Man-Mart} and does not rely on
the conjectural homotopy equivalence between $\Def_{\mathcal
F|_U}$ and $\MC_{\bullet}(\Eps nd^*(\Eps^\cdot)(U))$. More
precisely, we associate with any semicosimplicial DGLA ${\mathfrak
g}^\Delta$ a set-valued functor of Artin rings $Z^1_{\rm
sc}(\exp{\mathfrak g}^\Delta)$ together with an equivalence
relation $\sim$ on it, such that the quotient functor $ H^1_{\rm
sc}(\exp{\mathfrak g}^\Delta)=Z^1_{\rm sc}(\exp{\mathfrak
g}^\Delta)/\sim$ is an abstract version of $H^1_{\rm Ho}(X;\exp
\Eps nd^*(\Eps^\cdot))$. The latter is obtained, as a particular case, by
considering the \v{C}ech semicosimplicial Lie algebra ${\Eps nd^*(
\Eps^\cdot})(\mathcal U)$
\[
\xymatrix{ {\prod_i{\Eps nd^*(\Eps^\cdot})(U_i)}
\ar@<2pt>[r]\ar@<-2pt>[r] & {
\prod_{i<j}{\Eps nd^*(\Eps^\cdot})(U_{ij})}
      \ar@<4pt>[r] \ar[r] \ar@<-4pt>[r] &
      {\prod_{i<j<k}{\Eps nd^*(\Eps^\cdot})(U_{ijk})}
\ar@<6pt>[r] \ar@<2pt>[r] \ar@<-2pt>[r] \ar@<-6pt>[r]& \cdots}.
\]
Namely,
\[
H^1_{\rm Ho}(X;\exp \Eps nd^*(\Eps^\cdot))=\lim_{\substack{\longrightarrow \\ {\mathcal U}}}H^1_{\rm sc}(\exp{\Eps nd^*(\Eps^\cdot})(\mathcal U))
\]
and both sides coincide with $H^1_{\rm sc}(\exp \Eps
nd^*(\Eps^\cdot)(\U))$, for  an $\Eps nd^*(\Eps^\cdot)$-acyclic cover
of $X$. Next, we consider the Thom-Whitney model
$\Tot_{TW}{\mathfrak g}^\Delta$ for $\holim {\mathfrak g}^\Delta$
and show that there exists a commutative diagram of functors
\[ \xymatrix{
{\rm DGLA}^{\Delta_{\rm mon}}_{H^{\geq 0}}
\ar[rr]^{\Tot_{TW}}\ar[rd]_{H^1_{\rm sc}(\exp-)}
 &&{\rm DGLA}\ar[ld]^{\phantom{mi}\text{Maurer-Cartan}/\text{gauge}}\\
&\mathbf{Set}^{{\mathbf{Art}}_\K},& }
\]
where ${\rm DGLA}^{\Delta_{\rm mon}}_{H^{\geq 0}}$ is the category
of semicosimplicial DGLAs with no negative cohomology.  From the
point of view of $\infty$-groupoids, this can be seen as an
explicit description of the set $\pi_{\leq
0}(\MC_\bullet(\holim{\mathfrak g}^\Delta))$.
\bigskip

The paper is organized as follows: in Section 1 we dicuss deformations of coherent sheaves from a classical perspective and show how deformation data can be conveniently encoded into a \v{C}ech cohomology group with coefficient in a sheaf of DGLAs. In Section 2,  the functors $H^1_{\rm
sc}(\exp{\mathfrak g}^\Delta)$ and $H^1_{\rm Ho}(X;\exp {\mathcal L})$ are defined; next, in Sections 3 and 4, we recall the definition of the Thom-Whitney DGLA associated with $\g^\Delta$ and with its truncations $\g^{\Delta_{[m,n]}}$. Sections 5 and 6 are rather technical; namely Section 5 is devoted to a technical lemma on Maurer-Cartan elements in the Thom-Whitney  DGLAs $\Tot_{TW}(\g^{\Delta_{[0,1]}})$ and  $\Tot_{TW}(\g^{\Delta_{[0,2]}})$ and Section 6 to the proof of the isomorphism $H^1_{\rm
sc}(\exp{\mathfrak g}^{\Delta_{[0,1]}})$ and $\Def_{\Tot_{TW}(\g^{\Delta_{[0,1]}})}$. Finally, in Section 7, we are able 
 to prove our main result (Theorem \ref{th.iso funt}): under the cohomological hypotesis $H^{-1}(\g_2)=0$ there is a natural isomorphism of funtors
$\Def_{\Tot_{TW}(\g^{\Delta_{[0,2]}})} \cong H^1_{\rm sc}(\exp
\g^\Delta)$; moreover, if $H^{j}(\g_i)=0$ for all $i\geq 0$ and
$j<0$, then there is a natural isomorphism of functors
$\Def_{\Tot_{TW}(\g^{\Delta})} \cong H^1_{\rm sc}(\exp
\g^\Delta)$. In the concluding Section 8, we use this isomorphism to prove  that
infinitesimal deformations of a coherent sheaf ${\mathcal F}$ are controlled by the
DGLA of global sections of an acyclic resolution of $\Eps nd^*(\Eps^\cdot)$, where $\Eps^\cdot$ is a locally free resolution of $\mathcal F$.

\bigskip
While revising this paper, we became aware of \cite{yakutieli} where a similar construction is developed and investigated.
\bigskip

Throughout this paper we work on a fixed algebraically closed field $\K$ of 
characteristic zero; the symbol $\bf{Art}_{\K}$
 denotes the category of local Artinian $\K$-algebras $(A,{\mathfrak m}_A)$, with residue field $\K$.

\begin{acknowledgement}
We thank Marco Manetti for stimulating discussions on the subject 
and for useful comments and suggestions on the first version of this paper; d.i. thanks the Mathematical Department
\lq\lq Guido Castelnuovo\rq \rq , Sapienza Universit\`a di Roma for
the hospitality.

\end{acknowledgement}

\section{Infinitesimal deformations and sheaves of DGLAs} \label{sec.Example}
In this section, we study infinitesimal deformations of a coherent
sheaf $\mathcal F$  of $\O_X$-modules on a smooth projective
variety $X$ and explain
how these deformations can be naturally described in terms of a
sheaf of differential graded Lie algebras on $X$.

\smallskip

An infinitesimal deformation of the coherent sheaf
of $\O_X$-modules $\mathcal{F}$ over $A\in \bf{Art}_{\K}$ is given by a
coherent  sheaf $\mathcal{F}_A$ of $\O_X\otimes A$-modules on
$X\times \Spec A$, flat over $A$, with a morphism of sheaves $\pi:
\mathcal{F}_A \to \mathcal{F}$  inducing an isomorphism
$\mathcal{F}_A\otimes_A \K\cong \mathcal{F}$.

Two deformations $\mathcal{F}_A, \mathcal{F'}_A$ of the coherent
sheaf $\mathcal{F}$ over $A$ are isomorphic if there exists an
isomorphism of sheaves $f: \mathcal{F}_A \to \mathcal{F'}_A$, that
commutes with the morphisms to $\mathcal{F}$.
We denote by $\Def_{\mathcal{F}}: \bf{Art}_{\K}\to \bf{Set}$  the functor
of infinitesimal deformations of the sheaf $\mathcal{F}$.

\smallskip

We start by studying  infinitesimal deformations of a coherent
sheaf ${\mathcal F}$ of $\O_X$-modules on an affine variety $X$.
Let $X=\Spec R$, where $R$ is a Noetherian $\K$-algebra and let
$\mathcal F$ be the coherent sheaf associated with a finitely generated
$R$-module $M$; in this simple case, deformations of the sheaf
$\mathcal F$ reduce to deformations of the $R$-module $M$.

An infinitesimal deformation of the $R$-module $M$ over $A\in \bf{Art}_{\K}$ is given by a $R\otimes A$-module $M_A$, flat over $A$, with a morphism $\pi: M_A \to M$  inducing an isomorphism
$M_A\otimes_A \K\cong M$.
Two deformations $M_A$ and $ M'_A$ of the module $M$ over $A$ are isomorphic if there exists an
isomorphism of $R\otimes A$-modules $f: M_A \to M'_A$, that
commutes with the morphisms to $M$.

Next, let
\begin{equation} \label{presentazione}
\cdots \stackrel{d}{\longrightarrow }  R^{n_{1}}
\stackrel{d}{\longrightarrow }  R^{n_0}
\stackrel{d}{\longrightarrow } M \longrightarrow  0
\end{equation}
be a presentation of $M$ as $R$-module.
If $M_A$ is a deformation of $M$ over $A$, then it is an $A$-flat $R\otimes A$-module; therefore, flatness allows to lift relations between generators and to construct the exact sequence
\[\cdots \stackrel{d_A}{\longrightarrow }  R^{n_{1}}\otimes A
\stackrel{d_A}{\longrightarrow }  R^{n_0}\otimes A
\stackrel{d_A}{\longrightarrow } M_A \longrightarrow  0, \]
that reduces to (\ref{presentazione}) when tensored by $\K$ over $A$.
On the other hand, the datum of such an exact sequence assures flatness of the $R\otimes A$-module $M_A$ and so it defines a deformation of $M$ over $A$ (see \cite[par. 3]{Artin}, or \cite[Theorem A.31]{Sernesi} for details of these correspondences).
Moreover, if $M_A$ and $M'_A$ are isomorphic deformations of $M$ over $A$, the isomorphism between them lifts to an isomorphism between the correspondent deformed complexes and viceversa.

\smallskip

Next, we return to the global case of a coherent sheaf $\mathcal F$ of $\O_X$-modules on a smooth projective variety $X$.
Let
\[0 \longrightarrow \Eps^{-m} \stackrel{d}{\longrightarrow }
\cdots \stackrel{d}{\longrightarrow }  \Eps^{-1}
\stackrel{d}{\longrightarrow }  \Eps^0
\stackrel{d}{\longrightarrow } \mathcal F \longrightarrow  0 \]
be a global syzygy for ${\mathcal F}$, and  denote by $\Eps^{\cdot}$ the complex of
locally free sheaves
\[ (\Eps^{\cdot},d): \qquad\qquad  0\longrightarrow  \Eps^{-m}
 \stackrel{d}{\longrightarrow } \cdots\stackrel{d}{\longrightarrow }
  \Eps^{-1} \stackrel{d}{\longrightarrow }\Eps^0 \longrightarrow 0.\]
Let $\U=\{ U_i \}_{i\in I}$ be an affine\footnote{or Stein, if we work in the complex analytic category.} open cover of $X$, such
that every   sheaf of  $\Eps^{\cdot}$ is  free on  each $U_i$.

The Kodaira-Spencer approach to infinitesimal deformations of
$\mathcal F$ consists in deforming the sheaf $\mathcal F$ locally
in such a way that local deformations glue together to a global
sheaf, or equivalently, in view of the above discussion of the
affine case, in deforming the complex $(\Eps^{\cdot},d)$ on every
open set $U_i$ in such a way that these data glue together in
cohomology.

Following this approach, let us make explicit the deformation data:
the first datum is an element $l=\{l_i\}_i
\in \prod_i {\Eps nd}^1(\Eps^{\cdot})(U_i)\otimes \m_A$ defining, on every open set $U_i$, a complex
$(\Eps^{\cdot}|_{U_i}\otimes A,d + l_i)$ which is a deformation of the complex
$(\Eps^{\cdot}|_{U_i},d)$. Note that the condition for $(\Eps^{\cdot}|_{U_i}\otimes A,d + l_i)$ to be a complex
is the Maurer-Cartan equation:
\[dl_i+\frac{1}{2}[l_i,l_i]=0, \quad \mbox{for all } i\in I.\]
Also note that, by upper semicontinuity of cohomology, the complex $(\Eps^{\cdot}|_{U_i}\otimes A,d + l_i)$ is exact except possibly at
zero level.
To glue together the deformed local complexes $(\Eps^{\cdot}|_{U_i}\otimes A,d + l_i)$, we need to specify isomorphisms
between the deformed complexes on the double intersections of open
sets of the cover ${\mathcal U}$. Since these isomorphisms will have to be deformations of the identity, they will be of
the form
\[e^{m_{ij}}: (\Eps^{\cdot}|_{U_{ij}}\otimes A,d + l_j) \to
(\Eps^{\cdot}|_{U_{ij}}\otimes A,d + l_i),
\]
with
$m=\{m_{ij}\}_{i<j}\in \prod_{i<j} {\Eps
nd}^0(\Eps^{\cdot})(U_{ij})\otimes \m_A$. The compatibiliy with the differentials, i.e., the commutativity of the diagrams
\[ \xymatrix{
\Eps^{\cdot}|_{U_{ij}}\otimes A \ar[r]^{e^{m_{ij}}}\ar[d]_{d+l_j|_{U_{ij}}} &\Eps^{\cdot}|_{U_{ij}}\otimes A\ar[d]^{d+l_i|_{U_{ij}}}\\
\Eps^{\cdot}|_{U_{ij}}\otimes A \ar[r]^{e^{m_{ij}}} &\Eps^{\cdot}|_{U_{ij}}\otimes A
}
\]
can be written as
$d+l_i|_{U_{ij}}=e^{m_{ij}}(d+l_j|_{U_{ij}})e^{-m_{ij}}$, i.e., as
\[   l_i|_{U_{ij}}=e^{m_{ij}}*l_j|_{U_{ij}}, \quad \mbox{for all } i<j. \]
Finally, the above isomorphisms have to satisfy the cocycle condition up to homotopy.
Indeed, in order to obtain a deformation of ${\mathcal F}$, 
   we actually do not want to glue together the complexes $(\Eps^{\cdot}|_{U_i}\otimes A,d + l_i)$, but rather their cohomology sheaves. In other words, we require
$e^{m_{jk}} e^{-m_{ik}} e^{m_{ij}}$ to be
homotopic to the identity on triple intersections. Taking logarithm, what we require is that
$m_{jk} \bullet -m_{ik} \bullet m_{ij}$ is homotopy equivalent to
zero, i.e.,
\[
m_{jk} |_{U_{ijk}} \bullet - m_{ik}|_{U_{ijk}} \bullet m_{ij}|_{U_{ijk}} =
[d + l_j |_{U_{ijk}} , n_{ijk}],
\]
for some $n=\{n_{ijk}\}_{ijk} \in \prod_{i<j<k} {\Eps
nd}^{-1}(\Eps^{\cdot})(U_{ijk})$.
This homotopy cocycle equation is conveniently rewritten as
\[ m_{jk} |_{U_{ijk}} \bullet - m_{ik}|_{U_{ijk}} \bullet m_{ij}|_{U_{ijk}} =
d_{\Eps
nd^{*}(\Eps^{\cdot})}n_{ijk}+[l_j |_{U_{ijk}} , n_{ijk}].
\]

Next, let explain how the data introduced above are concretely
linked with deformations of the coherent sheaf $\mathcal F$ over
$A$. As the homotopy cocycle equation is satisfied, the local
$A$-flat sheaves of $\O_X|_{U_i}\otimes A$-modules ${\mathcal
F}_{A,U_i}:={\mathcal H}^*(\Eps^{\cdot}|_{U_i}\otimes A,d + l_i)$
glue together to give a global coherent sheaf ${\mathcal F}_A$
which is a deformation of ${\mathcal F}$. On the other hand, every
deformation ${\mathcal F}_A$ of the sheaf $\mathcal F$ can be
obtained in this way. Indeed, the resolution $(\Eps^{\cdot}, d)$
locally extends to projective resolutions
$(\Eps^{\cdot}\vert_{U_i}\otimes A, d+l_i)$ of ${\mathcal
F}_A\vert_{U_i}$; these deformed local resolutions are linked each
other on double intersections by isomorphisms of complexes lifting
the identity of ${\mathcal F}_A$ and the compositions of these
isomorphisms on triple intersections are homotopy to the identity,
since they lift the identity of ${\mathcal F}_A$ and liftings are
unique up to homotopy.

Let now ${\mathcal F}_A$ and ${\mathcal F'}_A$ be  isomorphic
deformations of the sheaf ${\mathcal F}$, associated with deformation data $(l,m)$
and $(l',m')$, respectively. The restriction to every open set
$U_i$ of the isomorphism between ${\mathcal F}_A$ and ${\mathcal F'}_A$ lifts to local isomorphisms
between the correspondent deformed  complexes. Since these isomorphisms
specialize to identities of $(\Eps^\cdot\vert_{U_i},d)$, they are of the form $e^{a_{i}}:
(\Eps^\cdot|_{U_i}\otimes A, d+l_i) \to (\Eps^\cdot|_{U_i}\otimes A, d+l'_i)$,
where $a=\{a_i\}_i \in \prod_i {\Eps nd}^0(\Eps^{\cdot})(U_i)\otimes
\m_A$. As above, compatibility with the differentials translates into the equations
\[e^{a_i} * l_i=l'_i, \quad \mbox{for all } i\in I. \]
Finally, since the local isomorphisms $e^{a_i}$ lift a global isomorphism in cohomology, the diagrams
\[ \xymatrix{
(\Eps^{\cdot}|_{U_{ij}}\otimes A,  d+l_j\vert_{U_{ij}})\ar[r]^{e^{m_{ij}}}\ar[d]_{e^{a_j}|_{U_{ij}}} &(\Eps^{\cdot}|_{U_{ij}}\otimes A,  d+l_i\vert_{U_{ij}})\ar[d]^{e^{a_i}|_{U_{ij}}}\\
(\Eps^{\cdot}|_{U_{ij}}\otimes A,  d+l'_j\vert_{U_{ij}}) \ar[r]^{e^{m'_{ij}}} &(\Eps^{\cdot}|_{U_{ij}}\otimes A,  d+l'_i\vert_{U_{ij}}),
}
\]
expressing compatibility with the gluing morphisms, commute in cohomology. Moreover, since the compositions
$e^{-m_{ij}} e^{-a_i} e^{ m'_{ij}} e^{a_j}$ lift the identity of ${\mathcal F}_A$ on double intersections and liftings are unique up to homotopy, these compositions are homotopy to identity and, reasoning as above, we find
\[ -m_{ij} \bullet -a_i|_{U_{ij}} \bullet m'_{ij} \bullet a_j
|_{U_{ij}} = d_{\Eps nd^{*}(\Eps^{\cdot})}b_{ij}+ [l_j|_{U_{ij}}, b_{ij}],\]
for some $b=\{b_{ij}\}_{i<j} \in \prod_{i<j} {\Eps nd}^{-1}(\Eps^{\cdot})(U_{ij})\otimes
\m_A$.
Viceversa, if for the deformation data $(l,m)$ and $(l',m')$ there exist $a=\{a_i\}_i \in \prod_i {\Eps nd}^0(\Eps^{\cdot})(U_i)\otimes \m_A$ and $b=\{b_{ij}\}_{i<j} \in \prod_{i<j} {\Eps nd}^{-1}(\Eps^{\cdot})(U_{ij})\otimes
\m_A$ that satisfy equations above, the local isomorphisms $e^{a_i}$ glue together in cohomology to give a global isomorphism of the correspondent deformed sheaves ${\mathcal F}_A$ and ${\mathcal F'}_A$.

Summing up, we have shown that in the Kodaira-Spencer approach,
infinitesimal deformations of the coherent sheaf ${\mathcal F}$
are controlled by the sheaf of DGLAs ${\Eps nd}^*(\Eps^{\cdot})$,
via the equations above. At the end of Section \ref{sec.main}, we
will apply techniques of semicosimplicial DGLAs developed in this
paper to recover the classical well known fact that the functor of
infinitesimal deformations of $\mathcal F$ has $\Ext^1({\mathcal
F},{\mathcal F})$ as tangent space and its obstructions are
contained in $\Ext^2({\mathcal F},{\mathcal F})$.
\begin{remark}
The above description of the functor of infinitesimal deformations of ${\mathcal F}$ is actually independent of the resolution chosen. Indeed, the DGLAs of the
endomorphisms of any two locally free resolutions of ${\mathcal F}$ are
quasi-isomorphic (see,e.g., \cite[Lemma 4.4]{Seidel.Thomas}).
\end{remark}
\begin{remark}
If the sheaf $\mathcal F$ is locally free, then we can take its
trivial resolution $0\to\mathcal F\to \mathcal F\to 0$; thus,  we
recover the well known fact that the infinitesimal deformations
of $\mathcal F$ are controlled by the sheaf $\Eps nd (\mathcal
F)$ of the endomorphism of
$\mathcal F$ , via the \v{C}ech functor $H^1(X,\Eps nd (\mathcal
F))$.
\end{remark}
\begin{remark}
Note that the results of this section actually hold under the hypotesis that ${\mathcal F}$ admists a global syzygy. This hypothesis is always satisfied, but in the general case the resolution is less obvious. Indeed, following Illusie \cite[Section 1.5]{Illusie}, for any sheaf $\mathcal F$ of ${\mathcal O}_X$-modules on a topological space $X$, one can construct the \emph{standard free resolution} of $\mathcal F$:
\[  \ldots \longrightarrow {\mathcal R}({\mathcal F})^{2}  \stackrel{D^2}{\longrightarrow} {\mathcal R}({\mathcal F})^{1}  \stackrel{D^1}{\longrightarrow}   {\mathcal R}({\mathcal F})^{0} \longrightarrow \mathcal F \longrightarrow 0.\]
Its terms are defined by recurrence: ${\mathcal R}({\mathcal F})^{0}$ is the free sheaf of $\O_X$-modules associated with the presheaf $U\mapsto \O_X(U)^{{\mathcal F}(U)}$, given on every open set $U\subset X$ by the free $\O_X(U)$-module generated by ${\mathcal F}(U)$; ${\mathcal R}({\mathcal F})^{j}$ is the free sheaf of $\O_X$-modules associated with the presheaf $U\mapsto \O_X(U)^{{\mathcal R}({\mathcal F})^{j-1}(U)}$, given on every open set $U\subset X$ by the free $\O_X(U)$-module generated by ${\mathcal R}({\mathcal F})^{j-1}(U)$.

To define morphisms $D^j$, let's write explicitly elements in ${\mathcal R}({\mathcal F})^{j}(U)$.
An element in ${\mathcal R}({\mathcal F})^{0}(U)$ is of the form $a^{i_0} \odot f_{i_0}$, where $a^{i_0}\in \O_X(U)$, $f_{i_0}\in \mathcal F (U)$, and we used the $\odot$ to denote the action of ${\mathcal O}_X(U)$ on the free $\O_X(U)$-module generated by $ \mathcal F (U)$, in order to distinguish it from the action of ${\mathcal O}_X(U)$ on the $\O_X(U)$-module $ \mathcal F (U)$.
Recursively,
an element in ${\mathcal R}({\mathcal F})^{j}(U)$ is of the form 
\[
a^{i_j}\odot a_{i_j}^{i_{j-1}} \odot\cdots \odot a^{i_0}_{i_1}\odot f_{i_0}
\]
where $a_{i_k}^{i_{k-1}}\in {\mathcal O}_X(U)$,  $f_{i_0}\in \mathcal F (U)$.
The differential of the resolution is defined as $D^j=\sum_{k=0}^{j}(-1)^i d^j_k$, where $d^j_k:{\mathcal R}({\mathcal F})^{j} \longrightarrow {\mathcal R}({\mathcal F})^{j-1}$ is defined by
\[
a^{i_j}\odot\cdots \odot a_{i_{k+1}}^{i_{k}}\odot a_{i_k}^{i_{k-1}}\odot \cdots \odot a^{i_0}_{i_1}\odot f_{i_0} 
\mapsto 
a^{i_j} \odot\cdots \odot a_{i_{k+1}}^{i_{k}} a_{i_k}^{i_{k-1}}\odot \cdots \odot a^{i_0}_{i_1}\odot f_{i_0} 
\]

The relevant fact is that the sequence of free sheaves of $\O_X$-modules $({\mathcal R}({\mathcal F})^{\cdot}, D^\cdot) \to \mathcal F$ is a resolution of $\mathcal F$ \cite[Theorem 1.5.3]{Illusie}. This construction can be done for every sheaf $\mathcal F$ of $\O_X$-modules on a topological space $X$; Illusie obtains it as an example of the even more general construction of the standard simplicial resolution of a pair of adjont functors  \cite[Section 1.5]{Illusie}.

\end{remark}

\section{Semicosimplicial DGLAs and the functor $H^1_{\rm sc}
(\exp{\mathfrak g}^\Delta)$} \label{sec.semicosimplicial dglas}

A \emph{semicosimplicial differential graded Lie algebra} is a
covariant functor $\mathbf{\Delta}_{\operatorname{mon}}\to
\mathbf{DGLA}$, from the category
$\mathbf{\Delta}_{\operatorname{mon}}$, whose objects are finite
ordinal sets and whose morphisms are order-preserving injective
maps between them, to the category of DGLAs. Equivalently, a
semicosimplicial DGLA ${\mathfrak g}^\Delta$ is a diagram
 \[
\xymatrix{ {{\mathfrak g}_0}
\ar@<2pt>[r]\ar@<-2pt>[r] & { {\mathfrak g}_1}
      \ar@<4pt>[r] \ar[r] \ar@<-4pt>[r] & { {\mathfrak g}_2}
\ar@<6pt>[r] \ar@<2pt>[r] \ar@<-2pt>[r] \ar@<-6pt>[r]&
\cdots},
\]
where each ${\mathfrak g}_i$ is a DGLA, and for each
$i>0$ there are $i+1$ morphisms of DGLAs
\[
\partial_{k,i}\colon {\mathfrak g}_{i-1}\to {\mathfrak
g}_{i},
\qquad k=0,\dots,i,
\]
such that
$\partial_{k+1,i+1}\partial_{l,i}=\partial_{l,i+1}\partial_{k,i}$,
for any
$k\geq l$.

A classical example is the following: given a sheaf ${\mathcal L}$ of DGLAs on a topological space $X$, and an open cover ${\mathcal U}$ of $X$, one has the \v{C}ech cosimplicial DGLA ${\mathcal L}({\mathcal U})$,
\[
\xymatrix{ {\prod_i\mathcal{L}(U_i)}
\ar@<2pt>[r]\ar@<-2pt>[r] & {
\prod_{i<j}\mathcal{L}(U_{ij})}
      \ar@<4pt>[r] \ar[r] \ar@<-4pt>[r] &
      {\prod_{i<j<k}\mathcal{L}(U_{ijk})}
\ar@<6pt>[r] \ar@<2pt>[r] \ar@<-2pt>[r] \ar@<-6pt>[r]& \cdots},
\]
where the morphisms $\partial_{k,i}$ are the restriction maps.

\begin{definition}
Let $\mathfrak g^{\Delta}$ be a semicosimplicial DGLA. The functor
\[Z^1_{\rm sc}(\exp \g^{\Delta})    :
\mathbf{Art}_{\mathbb K} \to \mathbf{Set}\] is defined, for all
$A\in \mathbf{Art}_{\mathbb K}$, by
\[ Z^1_{\rm sc}(\exp \g^{\Delta})(A)= \left\{ (l,m)\in
(\g_0^1 \oplus \g^0_1) \otimes \mathfrak m_A \left|
\begin{array}{l}
dl+\frac{1}{2}[l,l]=0,\\ \partial_{1,1}l=e^{m}*\partial_{0,1}l, \\
{\partial_{0,2}m} \bullet {-\partial_{1,2}m} \bullet
{\partial_{2,2}m} =dn+[\partial_{2,2}\partial_{0,1}l,n]\\
\qquad\qquad\qquad\qquad \text{for some $n\in {\mathfrak g}_2^{-1}\otimes{\mathfrak m}_A$}
\end{array} \right.\right\}.\]
\end{definition}

\begin{remark}\label{prop irrelevant stabilizer}
In DGLA theory, given a DGLA $L$ and a Maurer-Cartan element $x$ in $\MC_L (A)$, the set
\[ {\rm Stab}(x)= \{dh+[x,h] \mid  h \in L^{-1}\otimes
 {\mathfrak m}_A\}\]
 is called the \emph{irrelevant stabilizer} of $x$. Note that
 ${\rm Stab}(x) \subseteq {\rm stab}(x)$, where
${\rm stab}(x)= \{ a \in L^0\otimes \m_A \mid e^a*x=x \}$
is the stabilizer of $x$ under the gauge action of $L^0\otimes \m_A$ on $\MC_L (A)$.
Also note that,
for any $a \in L^0\otimes \m_A$,
$e^a e^{{\rm Stab}(x)}  e^{-a}= e^{{\rm Stab}(y)}$,  with
$y=e^a*x.$ \end{remark}

We now introduce an equivalence relation on the set $ Z^1_{\rm sc}
(\exp \g^{\Delta})(A)$ as follows: we say that two elements
$(l_0,m_0)$ and $(l_1,m_1) \in Z^1_{\rm sc} (\exp \g^{\Delta})(A)$
are equivalent under the relation  $\sim$ if  and only if there
exist elements $a \in \g^0_0\otimes \m_A$ and $b\in{\mathfrak
g}_1^{-1}\otimes{\mathfrak m}_A$ such that
\[
\begin{cases}
e^a * l_0=l_1\\
- m_0\bullet -\partial_{1,1}a \bullet m_1
\bullet \partial_{0,1}a=db+[\partial_{0,1}l_0,b].
\end{cases}
\]
\begin{remark}
The relation $\sim$ is actually an equivalence relation on
$Z^1_{\rm sc}(\exp \g^{\Delta})(A)$. First note that the set
$Z^1_{\rm sc}(\exp \g^{\Delta})(A)$ is closed under $\sim$.
Indeed, let $(l_0,m_0)$ and $(l_1,m_1) \in (\g_0^1\oplus \g_1^0)
\otimes \m_A$ be equivalent under $\sim$ via elements $a \in
\g^0_0\otimes \m_A$ and $b\in{\mathfrak g}_1^{-1}\otimes{\mathfrak
m}_A$, and  suppose that $(l_0,m_0) \in Z^1_{\rm sc} (\exp
\g^{\Delta})(A)$. Then $l_1=e^a*l_0$ satisfies the Maurer Cartan
equation and
\[ e^{m_1}*\partial_{0,1}l_1= e^{ \partial_{1,1}a \bullet m_0 \bullet (db+[\partial_{0,1}l_0,b])\bullet -\partial_{0,1}a}   e^{\partial_{0,1}a }
* \partial_{0,1} l_0=  e^{\partial_{1,1}a }
* \partial_{1,1} l_0 = \partial_{1,1}l_1. \]
Moreover,   an easy calculation, using relations between maps $\partial_{j,k}$ and Remark \ref{prop irrelevant stabilizer}, shows that
${\partial_{0,2}m_1} \bullet {-\partial_{1,2}m_1} \bullet {\partial_{2,2}m_1}$
is an element of the irrelevant stabilizer of $\partial_{2,2}\partial_{0,1}l_1$.

Secondly $\sim$ is an equivalent relation.
Reflexivity is trivial; for simmetry,
let $(l_0,m_0)$ and $(l_1,m_1)$ be equivalent via elements
$\ a \in \g^0_0\otimes \m_A$ and $b \in {\mathfrak g}_1^{-1}\otimes{\mathfrak m}_A$, then $e^{-a} * l_1=l_0$ and
$ - m_1\bullet \partial_{1,1}(a) \bullet m_0 \bullet-\partial_{0,1}(a)= \partial_{0,1}(a) \bullet -(db+[\partial_{0,1}l_0,b]) \bullet -\partial_{0,1}(a)$ is an element of the irrelevant stabilizer of $\partial_{0,1}l_1$, by Remark \ref{prop irrelevant stabilizer}.
Next, let $(l_0,m_0)\sim (l_1,m_1)$ via $a \in \g^0_0\otimes \m_A$ and $b\in{\mathfrak g}_1^{-1}\otimes{\mathfrak m}_A $, and $(l_1,m_1) \sim (l_2,m_2)$ via $\alpha  \in \g^0_0\otimes \m_A$ and $\beta\in{\mathfrak g}_1^{-1}\otimes{\mathfrak m}_A  $;
then, $e^{\alpha \bullet a}*l_0=l_2$ and
\[
- m_0\bullet  \partial_{1,1}(-(b\bullet a)) \bullet m_2 \bullet \partial_{0,1}(b\bullet a)  =  - m_0\bullet - \partial_{1,1} (a) \bullet - \partial_{1,1} (b) \bullet m_2 \bullet \partial_{0,1} (b)  \bullet \partial_{0,1} (a) =
\]
\[
- m_0\bullet - \partial_{1,1}(a) \bullet m_1\bullet (db+[\partial_{0,1}l_0,b])  \bullet \partial_{0,1}(a),
\]
by Remark \ref{prop irrelevant stabilizer}, it is an element of the irrelevant stabilizer of $\partial_{0,1}l_0$, therefore $\sim$ is transitive.
\end{remark}

\begin{definition}
Let $\g^\Delta$ be a semicosimplicial DGLA, the functor
\[ H^1_{\rm sc}(\exp \g^\Delta)  : \mathbf{Art}_{\mathbb K}
 \to \mathbf{Set}\]
is defined, for all $A\in \mathbf{Art}_{\mathbb K}$, by
\[H^1_{\rm sc}(\exp \g^\Delta)(A)= \frac{Z^1_{\rm sc}
(\exp \g^\Delta)(A)}{\sim}.\]
\end{definition}

\begin{remark}
Note that, if $\g^\Delta$ is a semicosimplicial  Lie
algebra, i.e., if all the DGLAs ${\mathfrak g}_i$ are concentrated in degree zero, then the functor
$H^1_{\rm sc}(\exp \g^\Delta)$ reduces to
the one defined in \cite{Fio-Man-Mart}.
\end{remark}

\begin{lemma}
The projection $ \pi: Z^1_{\rm sc} (\exp\g^{\Delta})
 \longrightarrow H^1_{\rm sc}(\exp \g^\Delta)$ is a
  smooth morphism of functors.
\end{lemma}
\begin{proof}
Let $ \beta: B  \longrightarrow A$ be a surjection in
$\Art_{\mathbb K} $, we prove that the map
$$
Z^1_{\rm sc} (\exp
\g^{\Delta}) (B) \longrightarrow H^1_{\rm sc}(\exp \g^\Delta) (B)\times_{H^1_{\rm sc}(\exp \g^\Delta)(A)}Z^1_{\rm sc} (\exp \g^{\Delta}) (A),
$$
induced by
\begin{center}
$\xymatrix{Z^1_{\rm sc} (\exp
\g^{\Delta})(B)  \ar[r]^\beta \ar[d]_\pi & Z^1_{\rm sc} (\exp
\g^{\Delta})(A)   \ar[d]^\pi \\
H^1_{\rm sc}(\exp \g^\Delta)(B) \ar[r]_\beta  & H^1_{\rm sc}(\exp \g^\Delta) (A),   \\ }$
\end{center}
is surjective.
Let $( [(l,m)],(l_0,m_0)) \in H^1_{\rm sc}(\exp \g^\Delta) (B)\times_{H^1_{\rm sc}(\exp \g^\Delta)(A)}Z^1_{\rm sc} (\exp \g^{\Delta}) (A)$, then
$(\beta l,\beta m)$ and $(l_0,m_0)$ are gauge equivalent in
$ Z^1_{\rm sc} (\exp\g^{\Delta})(A)$, i.e., there exist $a  \in \g^0_0\otimes \m_A$ such that
$
e^a * \beta l = l_0$ and $
- \beta m \bullet -\partial_{1,1}a \bullet m_0
\bullet \partial_{0,1}a=db+[\partial_{0,1}\beta l,b] $, for some $b\in{\mathfrak g}_1^{-1}\otimes \m_A$.
Let $\tilde{a} \in \g^0_0 \otimes m_B$ and $\tilde{b} \in{\mathfrak g}_1^{-1}\otimes \m_B$ be liftings of $a$ and $b$, respectively. The element $(e^{\tilde{a} }*l, \partial_{1,1} \tilde{a} \bullet m \bullet (d\tilde{b}+[\partial_{0,1}l,\tilde{b}]) \bullet -\partial_{0,1} \tilde{a}) \in
 Z^1_{\rm sc} (\exp \g^{\Delta})(B)$ is a pre-image of $( [(l,m)],(l_0,m_0))$.
 \end{proof}

Next, let $\mathcal L$  be a sheaf of DGLAs on a topological space $X$
and $\U=\{U_i \}_{i\in I}$ an open cover. Considering the \v{C}ech
cosimplicial DGLA ${\mathcal L}({\mathcal U})$, we can define the
functor $H^1_{\rm sc}(\exp{\mathcal L}(\U))$. This functor depends
on the cover $\U$, but as shown in the following Lemma, the limit
over open covers is a well defined functor:
\[ H_{\rm Ho}^1(X; \exp \mathcal{L})= \lim_{\substack{\longrightarrow \\ \U}}H^1_{\rm sc}(\exp\mathcal L(\U)): \Art \to \Set. \]

\begin{lemma}\label{lemma limite raffinamenti}
Let $\U=\{U_\alpha \}_{\alpha\in I}$ and $\U'=\{U'_\alpha
\}_{\alpha\in I'}$ be open covers of $X$ with $\U'$ refinement of
$\U$ and let $\phi, \psi: I'\to I$ two refinement maps. Then, the
induced morphisms $\rho_{\phi},
 \rho_{\psi}:H^1_{\rm sc}(\exp\mathcal L(\U)) \to H^1_{\rm sc} (\exp\mathcal L(\U'))$  coincide.
\end{lemma}

\begin{proof}
Both $\phi$ and $\psi$ induce, for all $A\in \mathbf{Art}_{\mathbb K}$, a
morphism $Z^1_{\rm sc}(\exp \mathcal L(\U))(A) \to Z^1_{\rm sc}(\exp
\mathcal L(\U'))(A)$, defined sending $(l_i,m_{ij})$ to
$\rho_{\phi}(l_i,m_{ij})= ({l_{\phi\alpha}}|_{U'_\alpha},
{m_{\phi\alpha, \phi \beta}}|_{U'_{\alpha \beta}})$ and $
\rho_{\psi}(l_i,m_{ij})= ({l_{\psi\alpha}}|_{U'_\alpha},
{m_{\psi\alpha, \psi \beta}}|_{U'_{\alpha \beta}})$,
respectively.
Therefore, it remains to prove that $\rho_{\phi}(l_i,m_{ij})  \sim
\rho_{\psi}(l_i,m_{ij}) $, for all $(l_i,m_{ij}) \in Z^1_{\rm
sc}(\exp \mathcal L(\U)) (A)$, i.e.,  for all $\alpha \in I'$, there
exists $a_\alpha \in \mathcal L^0 (U'_{\alpha})\otimes m_A$ such that
\[
\begin{cases}
e^{a_\alpha} * {l_{\phi\alpha}}|_{U'_\alpha}=
{l_{\psi\alpha}} |_{U'_\alpha} \\
- {m_{\phi\alpha, \phi \beta}}|_{U'_{\alpha \beta}} \bullet -
{a_\alpha} |_{U'_{\alpha \beta}}\bullet {m_{\psi\alpha, \psi
\beta}}|_{U'_{\alpha \beta}} \bullet {a_\beta}|_{U'_{\alpha
\beta}} \in {\rm Stab}(l_{\phi \beta}|_{U'_{\alpha \beta}}).
\end{cases}
\]
A simple computation shows that it is enough to choose
${a_\alpha}:= {m_{\psi\alpha, \phi \alpha}}|_{U'_{\alpha}}$, for all $\alpha$ in $I'$.
\end{proof}


\begin{remark}
Having introduced the limit $H^1_{\rm Ho}(\exp\mathcal L)$, for a
sheaf of DGLAs $\mathcal L$ on a topological space $X$, the
results of Section \ref{sec.Example} can be restated as follows:
the functor of infinitesimal deformations of a coherent sheaf
${\mathcal F}$ on a projective manifold $X$ is
\[ \Def_{\mathcal F}\cong H^1_{\rm Ho}(X;\exp {\Eps}nd^*(\Eps^\cdot)),
\]
where $\Eps^\cdot$ is a locally free resolution of ${\mathcal F}$.

\end{remark}

 The example
of coherent sheaves on projective manifolds together with the DGLA
approach to deformation theory suggests that the functors of Artin
rings $H^1_{\rm sc}(\exp \g^\Delta)$ could actually be isomorphic
to functors $\Def_{L(\g^\Delta)}$ for some DGLA $L(\g^\Delta)$
canonically associated with $\g^\Delta$. We are going to show
that, under the cohomological hypothesis $H^{-1}(\g_2)=0$, it is
indeed so. More precisely, we are going to prove that, if
$H^{-1}(\g_2)=0$, then the functor  of Artin rings $H^1_{\rm
sc}(\exp \g^\Delta)$ is isomorphic to the deformation functor
associated with the Thom-Whitney DGLA of the truncation
${\mathfrak g}^{\Delta_{[0,2]}}$.

\section{The Thom-Whitney DGLA $\Tot_{TW}({\mathfrak g}^\Delta)$} \label{sec.def functors}
Let ${\mathfrak g}^\Delta$ be a semicosimplicial DGLA. The maps
\[
\partial_i=\partial_{0,i}-\partial_{1,i}+\cdots+(-1)^{i}
\partial_{i,i}
\]
endow the vector space $\bigoplus_i{\mathfrak g}_i$ with the
structure of a differential complex. Moreover, being a
DGLA, each ${\mathfrak g}_i$ is in particular a differential
complex
\[
{\mathfrak g}_i=\bigoplus_j {\mathfrak g}_i^j; \qquad
d_i\colon {\mathfrak g}_i^j\to{\mathfrak g}_i^{j+1}
\]
and since the maps $\partial_{k,i}$ are morphisms of DGLAs,
the space
\[
{\mathfrak g}^\bullet_\bullet=\bigoplus_{i,j}{\mathfrak
g}_i^j
\]
has a natural bicomplex structure. The associated total complex
\[({\rm Tot}({\mathfrak g}^\Delta),d_{\Tot})\quad\text{where}\quad
{\rm Tot}({\mathfrak g}^\Delta)=\bigoplus_{i}{\mathfrak
g}_i[-i],\quad d_{\Tot}=\sum_{i,j}\partial_i+(-1)^jd_j\] has no
natural DGLA structure. Yet there is an other bicomplex  naturally
associated with a semicosimplicial DGLA, whose total complex is naturally a
DGLA.

For every $n\ge 0$, denote by $\Omega_n$ the differential graded
commutative algebra of polynomial differential forms on the
standard $n$-simplex $\Delta^n$:
\[ \Omega_n=\frac{{\mathbb K}[t_0,\ldots,t_n,dt_0,
\ldots,dt_n]}{(\sum t_i-1,\sum dt_i)}.\]
Denote
by $\delta^{k,n}\colon \Omega_n\to \Omega_{n-1}$, $k=0,\ldots,n$,
the face maps; then, one has  natural morphisms of bigraded DGLAs
\[ \delta^{k,n}\colon \Omega_n\otimes \mathfrak{g}_n\to
\Omega_{n-1}\otimes\mathfrak{g}_n,\qquad
\partial_{k,n}\colon \Omega_{n-1}\otimes \mathfrak{g}_{n-1}\to
\Omega_{n-1}\otimes\mathfrak{g}_{n},\]
for every $0\le k\le n$.

The Thom-Whitney bicomplex
is defined as
\[
C^{i,j}_{TW}(\mathfrak{g}^\Delta)
=\{ (x_n)_{n\in {\mathbb
N}}\in \bigoplus_n \Omega_n^i\otimes {\mathfrak g}_n^j
\mid \delta^{k,n}x_n=
\partial_{k,n}x_{n-1}\quad \forall\; 0\le k\le n\},
\]
where  $\Omega_n^i$ denotes the degree $i$ component of $\Omega_n$.
Its total complex is a DGLA, called the  \emph{Thom-Whitney DGLA},
and it is denoted by
$\operatorname{Tot}_{TW}(\mathfrak{g}^\Delta)$; denote by
$d_{TW}$ the differential of the Thom-Whitney DGLA.
It is a remarkable fact that the integration maps
\[ \int_{\Delta^n}\otimes \operatorname{Id}\colon
\Omega_n\otimes \mathfrak{g}_n\to
\K[n]\otimes \mathfrak{g}_n=\mathfrak{g}_n[n]\]
give a quasi-isomorphism of differential complexes
\[
I\colon (\operatorname{Tot}_{TW}({\mathfrak
g}^\Delta), d_{TW})\to (\operatorname{Tot}({\mathfrak
g}^\Delta),d_{\rm Tot}).
\]
Moreover, Dupont has described in \cite{dupont1,dupont2} an
explicit morphism of differential complexes
\[
E\colon \operatorname{Tot}({\mathfrak
g}^\Delta)\to \operatorname{Tot}_{TW}({\mathfrak
g}^\Delta)
\]
and an explicit homotopy \[
h\colon \operatorname{Tot}_{TW}({\mathfrak
g}^\Delta)\to \operatorname{Tot}_{TW}({\mathfrak
g}^\Delta)[-1]
\]
such that
\[
IE={\rm Id}_{{\rm Tot}({\mathfrak
g}^\Delta)};\qquad
EI-{\rm Id}_{{\rm Tot}_{TW}({\mathfrak
g}^\Delta)}=[h,d_{TW}].
\]
We also refer to the papers \cite{chenggetzler,getzler,navarro}
for the explicit description of $E,h$ and for the proof of the
above identities. Here, we point out that $E$ and $h$ are defined
in  terms of integration over standard simplexes and
multiplication with canonical differential forms: in particular,
the construction of $\operatorname{Tot}_{TW}({\mathfrak
g}^\Delta)$, $\operatorname{Tot}({\mathfrak g}^\Delta)$, $I$, $E$
and $h$ is functorial in the category
$\mathbf{DGLA}^{\Delta_{\operatorname{mon}}}$ of
semicosimplicial DGLAs.\\

Recall that with a DGLA $L$ there is
a canonically associated deformation functor $\Def_L$, defined as the
solutions of Maurer-Cartan equation modulo gauge action (or, equivalently, modulo homotopy equivalence). Moreover, the tangent space to $\Def_L$ is $H^1(L)$ and obstructions live in $H^2(L)$. Thus,
with a semicosimplicial DGLA $\g^{\Delta}$ is also associated the
deformation functor $\Def_{{\Tot}_{TW}(\g^{\Delta})}$; its tangent space is
\[
T\Def_{{\Tot}_{TW}(\g^{\Delta})}\cong H^1({\Tot}_{TW}(\g^{\Delta}))\cong H^1({\Tot}(\g^{\Delta}))
\]
and obstructions live in
\[
H^2({\Tot}_{TW}(\g^{\Delta}))\cong H^2({\Tot}(\g^{\Delta})).
\]
\smallskip

Let $\mathbf{\Delta}^+_{\operatorname{mon}} $ the
category obtained by adding the empty set $\emptyset$ to the
category $\mathbf{\Delta}_{\operatorname{mon}}$. An
\emph{augmented semicosimplicial differential graded Lie algebra}
is a covariant functor $\mathbf{\Delta}^+_{\operatorname{mon}}\to
\mathbf{DGLA}$, from the category
$\mathbf{\Delta}^+_{\operatorname{mon}}$ to the category of DGLAs.
Equivalently, an augmented semicosimplicial DGLA ${\mathfrak
g}^{\Delta^+}$ is a diagram
\[
\xymatrix{{{\mathfrak g}_{-1}}\ar[r] &{{\mathfrak g}_0}
\ar@<2pt>[r]\ar@<-2pt>[r] & { {\mathfrak g}_1}
      \ar@<4pt>[r] \ar[r] \ar@<-4pt>[r] & { {\mathfrak g}_2}
\ar@<6pt>[r] \ar@<2pt>[r] \ar@<-2pt>[r] \ar@<-6pt>[r]& \cdots},
\]
where the truncated diagram ${\g}^\Delta$
\[
\xymatrix{{{\mathfrak g}_0} \ar@<2pt>[r]\ar@<-2pt>[r] & {
{\mathfrak g}_1}
      \ar@<4pt>[r] \ar[r] \ar@<-4pt>[r] & { {\mathfrak g}_2}
\ar@<6pt>[r] \ar@<2pt>[r] \ar@<-2pt>[r] \ar@<-6pt>[r]& \cdots}
\]
is a semicosimplicial DGLA and
\[
\partial_{0,0}\colon \g_{-1}\to \g_{0}
\]
is a DGLA morphism such that
$\partial_{0,1}\partial_{0,0}=\partial_{1,1}\partial_{0,0}$.

\begin{remark}\label{rm.aum}
There is a morphism of DGLAs
\begin{align*}
{\mathfrak g}_{-1} &\to \Tot_{TW}({\mathfrak g}^\Delta)\\
              x &\mapsto(\partial_{0,0}x,\ 
\partial_{1,1}\partial_{0,0}x,\ 
\partial_{2,2}\partial_{1,1}\partial_{0,0}x,\dots);
\end{align*}
the image of $x$ is an element in
$\Tot_{TW}({\mathfrak g}^\Delta)$ because of equations
$\partial_{1,1}\partial_{0,0}=\partial_{0,1}\partial_{0,0}$  and
$\partial_{k+1,i+1}\partial_{l,i}=\partial_{l,i+1}\partial_{k,i}$,
for any $k\geq l$. This morphism is obtained as the composition of
the natural inclusion $ {\mathfrak g}_{-1} \hookrightarrow
\Tot({\mathfrak g}^\Delta)$ with the morphism $E: \Tot({\mathfrak
g}^\Delta) \to  \Tot_{TW}({\mathfrak g}^\Delta)$. The existence of the DGLA
 morphism ${\mathfrak g}_{-1}\to \Tot_{TW}({\mathfrak g}^\Delta)$
 is not surprising; indeed, it is induced by the natural morphism
 $\lim{\mathfrak g}^\Delta\to \holim{\mathfrak g}^\Delta$.
\end{remark}

We use augmentation to link the Thom-Whitney DGLA of the \v{C}ech semicosimplicial DGLA of a sheaf of DGLAs with the DGLA of global sections of an acyclic resolution of the sheaf. This result is a translation of Theorem 7.2 in \cite{Fio-Man-Mart} in terms of the Thom-Whitney DGLA.

We recall that if ${\mathcal L}$ is a sheaf of DGLAs on a topological
space $X$ and ${\mathcal U}$ is an open cover of $X$, the associated \v{C}ech
semicosimplicial differential graded Lie algebra is:
\[ \mathcal{L}(\mathcal{U}):\quad
\xymatrix{ {\prod_i\mathcal{L}(U_i)}
\ar@<2pt>[r]\ar@<-2pt>[r] & {
\prod_{i<j}\mathcal{L}(U_{ij})}
      \ar@<4pt>[r] \ar[r] \ar@<-4pt>[r] &
      {\prod_{i<j<k}\mathcal{L}(U_{ijk})}
\ar@<6pt>[r] \ar@<2pt>[r] \ar@<-2pt>[r] \ar@<-6pt>[r]& \cdots}.
\]
A morphism
$\varphi\colon{\mathcal L}\to {\mathcal A}$ of sheaves of DGLAs is
a quasi-isomorphism if it is a quasi-isomorphism of sheaves of
differential complexes, i.e., if it induces linear isomorphisms
between the cohomology sheaves,
\[{\mathcal H}^*(\varphi)\colon{\mathcal H}^*({\mathcal
L})\xrightarrow{\sim} {\mathcal H}^*({\mathcal
A}).\]
Moreover, if ${\mathcal A}^k$ is an acyclic sheaf for any $k$, then $\varphi\colon{\mathcal L}\to {\mathcal A} $ is called an acyclic resolution of ${\mathcal L}$.

\begin{theorem} \label{th.semicosimpl vs dgla}
Let $X$ be a paracompact Hausdorff topological space, $\mathcal L$
a sheaf of differential graded Lie algebras on $X$, and
$\varphi\colon {\mathcal L}\to {\mathcal A}$ an acyclic
resolution. Also let $A={\mathcal A}(X)$ be the DGLA of global
sections of ${\mathcal A}$. Then, if $\mathcal U$ is an open
cover of $X$ which is acyclic with respect to both ${\mathcal
L}$ and ${\mathcal A}$, the DGLA
${\Tot_{TW}}({\mathcal L}(\mathcal{U}))$ is naturally
quasi-isomorphic to the DGLA $A$.
\end{theorem}

\begin{proof}
The natural inclusion $A\to \mathcal{A}(\mathcal{U})$ gives an
augmented semicosimplicial DGLA, and so it induces a morphism of DGLAs $A
\to \Tot_{TW}(\mathcal{A}(\mathcal{U}))$, that is the composition of the natural inclusion $A\to \Tot({\mathcal A}({\mathcal U}))$ with the quasi-isomorphism $E:\Tot({\mathcal A}({\mathcal U}))\to\Tot_{TW}(\mathcal{A}(\mathcal{U}))$, by Remark~\ref{rm.aum}.
Since the sheaves ${\mathcal A}^k$ are acyclic and
${\mathcal U}$-acyclic, and $A^k=H^0(X;\mathcal A^k)$, the inclusion $A \to
\Tot(\mathcal{A}(\mathcal{U}))$ is a
quasiisomorphism. Indeed, we have a natural identification
$H^*({\operatorname{Tot}}(\mathcal{A}(\mathcal{U})))={\mathbb
H}^*(X;{\mathcal A})$, and the spectral sequence abutting to the
hypercohomology of $X$ with coefficients in ${\mathcal A}$
degenerates at $E_2$, giving
\[
{\mathbb H}^k(X;{\mathcal
A})=\bigoplus_{p+q=k}E_2^{p,q}=E_2^{k,0}=
H^k(A).
\]
Then,  $A
\to \Tot_{TW}(\mathcal{A}(\mathcal{U}))$ is a quasi-isomorphism of DGLAs.

The morphism $\varphi\colon {\mathcal L}\to{\mathcal A}$ induces a
morphism of semicosimplicial DGLAs
\[ \varphi\colon \mathcal{L}(\mathcal{U})\to
\mathcal{A}(\mathcal{U}), \] and a morphism of complexes
\[ \varphi\colon \Tot_{TW}(\mathcal{L}(\mathcal{U})) \to
\Tot_{TW}(\mathcal{A}(\mathcal{U})).\]
Since the open
cover ${\mathcal U}$ is ${\mathcal L}$-acyclic, the cohomology
of the total complex $\Tot(\mathcal{L}(\mathcal{U}))$ is naturally
identified with the hypercohomology of $X$ with coefficients in
${\mathcal L}$,
\[
H^*(\Tot(\mathcal{L}(\mathcal{U})))\cong {\mathbb
H}^*(X;{\mathcal L}),
\]
and the induced linear map
\[
H^*(\varphi)\colon
H^*(\Tot(\mathcal{L}(\mathcal{U})))\to
H^*(\Tot(\mathcal{A}(\mathcal{U})))
\]
is identified with the linear map
\[
{\mathbb H}^*(\varphi)\colon {\mathbb H}^*(X;{\mathcal
L})\rightarrow{\mathbb H}^*(X;{\mathcal
A})
\]
induced in hypercohomology. Since, by hypothesis, $\varphi$ is a
quasi-isomorphism of sheaves of DGLAs, the induced map in
hypercohomology is an isomorphism, and so the morphism
$\varphi\colon \Tot(\mathcal{L}(\mathcal{U})) \to \Tot(\mathcal{A}(\mathcal{U}))$ is a
quasi-isomorphism of complexes.

Via the composition with quasi-isomorphisms $E$ and $I$ between the total complex and the Thom-Whitney total complex of a semicosimplicial DGLA, the morphism $\varphi$ induces a quasi-isomorphism of DGLAs
\[ \Tot_{TW}(\mathcal{L}(\mathcal{U})) \to
\Tot_{TW}(\mathcal{A}(\mathcal{U})).\]

Therefore, we have the
chain of quasi-isomorphisms of DGLAs
\[\Tot_{TW}(\mathcal{L}(\mathcal{U}))\xrightarrow{\sim}
\Tot_{TW}(\mathcal{A}(\mathcal{U}))
\xleftarrow{\sim}A.
\]
\end{proof}

\section{Truncations}
Let
\[ \mathfrak{g}^{\Delta}: \ \  \xymatrix{ \mathfrak{g}_0
\ar@<2pt>[r]\ar@<-2pt>[r] & \mathfrak{g}_1
      \ar@<4pt>[r] \ar[r] \ar@<-4pt>[r] & \mathfrak{g}_2
\ar@<6pt>[r] \ar@<2pt>[r] \ar@<-2pt>[r] \ar@<-6pt>[r] &
\mathfrak{g}_3 \ar@<8pt>[r] \ar@<4pt>[r] \ar[r]
\ar@<-4pt>[r] \ar@<-8pt>[r] &
\ldots }
\]
be a semicosimplicial DGLA. Let $m_1 \in \N$ and $m_2 \in \N \cup
\{\infty\}$ with $m_1 \leq m_2$, we denote by
$\g^{\Delta_{[m_1,m_2]}}$ the \emph{truncated between levels $m_1$
and $m_2$} semicosimplicial DGLA defined by
\[
(\g^{\Delta_{[m_1,m_2]}})_n=
\begin{cases}
\g_n &\text{for }m_1 \leq n\leq m_2\\
0&\text{otherwise},
\end{cases}
\]
with the obvious maps
$\partial_{k,i}^{[m_1,m_2]}=\partial_{k,i}$, for $m_1 < i\leq m_2$,
and $\partial_{k,i}^{[m_1,m_2]}=0$, otherwise. For any positive
integers $m_1,m_2,r_1,r_2$, such that $r_i\leq m_i$, the map
$\Id_{[m_1,r_2]}\colon \g^{\Delta_{[m_1,m_2]}}\to
\g^{\Delta_{[r_1,r_2]}}$ given by
\[
\Id_{[m_1,r_2]}\biggr\vert_{(\g^{\Delta_{[m_1,m_2]}})_n}=
\begin{cases}
\Id_{\g_n}&\text{if } m_1 \leq n \leq r_2\\
0&\text{otherwise}.
\end{cases}
\]
is a morphism of semicosimplicial DGLAs; it induces the  natural
morphism of complexes $\phi: \Tot(\g^{\Delta_{[m_1,m_2]}})\to
\Tot(\g^{\Delta_{[r_1,r_2]}})$ and the natural morphism of DGLAs
$\psi: \Tot_{TW}(\g^{\Delta_{[m_1,m_2]}})\to
\Tot_{TW}(\g^{\Delta_{[r_1,r_2]}})$. Note that we have an homotopy
commutative diagram of complexes
\[ \xymatrix{ \Tot(\g^{\Delta_{[m_1,m_2]}}) \ar@<1ex>[d]^{E}
\ar[r]^{\phi}& \Tot(\g^{\Delta_{[r_1,r_2]}}) \ar@<1ex>[d]^{E} \\
\Tot_{TW}(\g^{\Delta_{[m_1,m_2]}}) \ar[r]^{\psi}
 \ar@<1ex>[u]^{I}& \Tot_{TW}(\g^{\Delta_{[r_1,r_2]}})
 \ar@<1ex>[u]^{I}.}   \]

\begin{proposition} \label{prop.tre suff}
Let $\g^{\Delta}$ be a semicosimplicial DGLA such that
$H^j(\g_i)=0$, for all $i\geq 0$ and $j<0$. Then, the morphism
$\Id_{[0,2]}$ induces a natural isomorphism of functors:
$$
\Def_{\Tot_{TW}(\g^\Delta) } \xrightarrow{\sim}
\Def_{\Tot_{TW}(\g^{\Delta_{[0,2]}})}.
 $$
\end{proposition}

\begin{proof}
It is a well known fact (see, e.g., \cite{Man Pisa}  for a
proof), that a DGLA morphism which is surjective on $H^0$, bijective on $H^1$ and injective on $H^2$ induces an
isomorphism between the associated deformation functors. Since the above
homotopy commutative diagram identifies $H^*(\psi)$ with
$H^*(\phi)$,  it is enough to prove that $H^0(\phi)$ is
surjective, $H^1(\phi)$ is bijective and $H^2(\phi)$ is injective. This
is easily checked by looking at the spectral sequences
associated with double complexes of ${\mathfrak g}^\Delta$ and ${\mathfrak g}^{\Delta_{[0,2]}}$.
\end{proof}

\begin{remark} \label{rmk. truncation}
Observe that, for any semicosimplicial DGLA $\g^{\Delta}$, we have
$Z^1_{\rm sc}(\exp \g^\Delta)= Z_{\rm sc}^1 (\exp
\g^{\Delta_{[0,2]}})$ and $H^1_{\rm sc}(\exp \g^\Delta)= H^1_{\rm
sc}(\exp \g^{\Delta_{[0,2]}})$. Moreover, the inclusion     $Z^1_{\rm
sc}(\exp \g^\Delta)\hookrightarrow Z^1_{\rm sc} (\exp
\g^{\Delta_{[0,1]}})$ induces an injective
 map $H^1_{\rm sc}(\exp \g^\Delta) \hookrightarrow H^1_{\rm sc}
 (\exp \g^{\Delta_{[0,1]}})$.
\end{remark}
\begin{remark} \label{rmk. defN banale}
For later use, we point out that, if $\g^{\Delta}$ is a semicosimplicial DGLA with $H^{-1}(\g_2)=0$, then
$\Def_{\Tot_{TW}(\g^{\Delta_{[2,2]}})}$ is trivial. Indeed, $H^{1}(
\Tot_{TW}( \g^{ \Delta_{[2,2]}}) )= H^{-1}(\g_2)=0$.
\end{remark}
\begin{remark}Note that, by the definition of $H^1_{\rm sc}(\exp\g^\Delta)$ it follows that, if $H^{-1}(\g_2)=0$, then
\[
TH^1_{\rm sc}(\exp\g^\Delta)=H^1({\Tot}(\g^{\Delta_{[0,2]}})).
\]
Hence, the two functors of Artin rings $H^1_{\rm sc}(\exp\g^\Delta)$ and $\Def_{{\Tot}_{TW}(\g^{\Delta_{[0,2]}})}$ have naturally isomorphic tangent spaces when $H^{-1}(\g_2)=0$. We will show in Section \ref{sec.main} that in this case these two functors are actually isomorphic.
\end{remark}

\section{A lemma on Maurer-Cartan elements}
We will now give an explicit description of the solutions of Maurer-Cartan equation for the DGLAs
$\Tot_{TW}(\g^{\Delta_{[0,1]}})$ and $\Tot_{TW}(\g^{\Delta_{[0,2]}})$. Our main tool will be the following general result \cite[Proposition 7.2]{Fio-Man}:
\begin{lemma}\label{lemma FIO-MAN decomposizione}
Let $(L,d,[~,~])$ be a differential graded Lie algebra such that:
\begin{enumerate}

\item $L=M\oplus C\oplus D$ as graded vector spaces.

\item $M$ is a differential graded subalgebra of $L$.

\item  $d\colon C\to D[1]$ is an isomorphism of graded vector
spaces.
\end{enumerate}

Then, for every $A\in \mathbf{Art}_{\K}$ there exists a bijection
\[ \alpha\colon \MC_M(A)\times (C^0\otimes
\mathfrak{m}_A)\mapor{\sim}\MC_L(A),\qquad (x,c)\mapsto e^c\ast
x.\]
\end{lemma}
As almost immediate corollaries we obtain:
\begin{proposition}\label{prop.mc[01]}
Let $\g^{\Delta}$ be a semicosimplicial DGLA. Then, for every $A\in \mathbf{Art}_{\K}$, the solutions of the Maurer-Cartan equation for the Thom-Whitney DGLA
$\Tot_{TW}(\g^{\Delta_{[0,1]}})\otimes \m_A$ are of the form
$(x, e^{p(t)}* \partial_{0,1}x)$,
where $x \in \MC_{\g_0}(A)$ and  $p(t) \in (\g_1^0[t]\cdot t) \otimes \m_A$. The elements $x, p$ are uniquely determined, and they satisfy
\begin{equation}\label{eq.face conditions01}
\partial_{1,1}x=e^{p(1)}*\partial_{0,1}x.
\end{equation}
\end{proposition}

\begin{proof}
Notice that $\Tot_{TW}(\g^{\Delta_{[0,1]}})$ is a sub-DGLA of
$\g_0\oplus\Omega_1\otimes\g_1$.  Then, Lemma \ref{lemma FIO-MAN decomposizione}
with the decomposition of $\Omega_1\otimes \g_1$ given by
\[ M=\g_1, \qquad C= \g_1[t]\cdot t, \qquad D=dC \]
 tells us that every solution of the
Maurer-Cartan equation for $\Tot_{TW}(\g^{\Delta_{[0,1]}})\otimes
\m_A$ is of the form specified above.
\end{proof}

\begin{proposition}
Let $\g^{\Delta}$ be a semicosimplicial DGLA. Then, for every
$A\in \mathbf{Art}_{\K}$,  the solutions of the Maurer-Cartan
equation for the Thom-Whitney DGLA
$\Tot_{TW}(\g^{\Delta_{[0,2]}})\otimes \m_A$ are of the form
\[
 (x, e^{p(t)}* \partial_{0,1}x, e^{q(s_0,s_1)+ r(s_0,s_1,ds_0,ds_1)}*\partial_{0,2}\partial_{0,1}x),
\]
where $x \in \MC_{\g_0}(A)$,  $p(t) \in (\g_1^0[t]\cdot t)\otimes \m_A$, $q(s_0,s_1) \in (\g_2^0[s_0,s_1]\cdot s_0+ \g_2^0[s_0,s_1]\cdot s_1)\otimes \m_A$ and $r(s_0,s_1,ds_0,ds_1) \in (\g_2^{-1}[s_0,s_1]\cdot
s_0ds_1) \otimes \m_A$. The elements $x,p,q,r$ are uniquely determined, and they satisfy
\begin{equation}\label{eq.face conditions}
\begin{cases}
\partial_{1,1}x=e^{p(1)}*\partial_{0,1}x, \\
\partial_{0,2}p(t)= q(0,t), \\
 \partial_{1,2}p(t)= q(t,0),   \\
e^{(- \partial_{2,2}p(t))\bullet(q(t,1-t)+ r(t,1-t,dt))\bullet(-q(0,1))}*\partial_{2,2}\partial_{0,1}x=\partial_{2,2}\partial_{0,1}x.
\end{cases}
\end{equation}
\end{proposition}

\begin{proof}
Since $\Tot_{TW}(\g^{\Delta_{[0,2]}})$ is a sub-DGLA of
$\g_0\oplus\Omega_1\otimes\g_1\oplus\Omega_2\otimes\g_2$,
applying Lemma \ref{lemma FIO-MAN decomposizione}
with the decomposition of $\Omega_2\otimes \g_2$ given by
\[ M=\g_2,\qquad C=\g_2[s_0,s_1]\cdot s_0+ \g_2[s_0,s_1]\cdot
s_1 +\g_2[s_0,s_1]\cdot s_0ds_1 ,\qquad D=dC  \]
we obtain  that every solution of the
Maurer-Cartan equation for $\Tot_{TW}(\g^{\Delta_{[0,2]}})\otimes \m_A$ is of
the form
\[(x, e^{p(t)}* y, e^{q(s_0,s_1)+ r(s_0,s_1,ds_0,ds_1)}*z),\]
with the face conditions
\[
y=\partial_{0,1}x;\qquad z=\partial_{0,2}\partial_{0,1}x.
\]
The first relations in (\ref{eq.face conditions}) are a direct consequence
of face conditions and uniqueness. The last one is obtained as
follows. The last face condition is
\[ \partial_{2,2}(e^{p(t)}* \partial_{0,1}x) = e^{q(t,1-t)+ r(t,1-t,dt)}*\partial_{0,2}\partial_{0,1}x;\]
using  the other face conditions and  relations
 between maps $\partial_{k,i}$, we  obtain that
\[\partial_{2,2} \partial_{0,1}x =
\partial_{0,2}\partial_{1,1} x = \partial_{0,2}(e^{p(1)}* \partial_{0,1}x )= e^{q(0,1)}*\partial_{0,2}\partial_{0,1}x.\] Then, the above equation
becomes
\[ e^{\partial_{2,2}p(t)}* \partial_{2,2}\partial_{0,1}x =
 e^{(q(t,1-t)+ r(t,1-t,dt))\bullet (-q(0,1))}* \partial_{2,2}\partial_{0,1}x. \]
 \end{proof}

\section{The isomorphism $H^1_{\rm sc}(\exp \g^{\Delta_{[0,1]}}) \cong
\Def_{ \Tot_{TW}(\g^{ \Delta_{[0,1]} })}$}

\begin{proposition} \label{prop.mappa phi1}
Let $\g^\Delta$ be a semicosimplicial DGLA. The map
\[
\Phi_{[0,1]}:\MC_{ \Tot_{TW}(\g^{\Delta_{[0,1]}}) }(A)
\to (\g_0^1\oplus \g_1^0  ) \otimes \mathfrak m_A,
\]
given by
\[
(x, e^{p(t)}* \partial_{0,1}x)
\mapsto (x, p(1) ),
\]
induces a natural transformation of functors of Artin rings
\[ \Def_{\Tot_{TW}(\mathfrak{g}^{\Delta_{[0,1]}})}
\to H^1_{\rm sc}(\exp \g^{\Delta_{[0,1]}}).\]
\end{proposition}

\begin{proof}
Clearly, if $(x, e^{p(t)}* \partial_{0,1}x) \in
\MC_{ \Tot_{TW}(\g^{\Delta_{[0,1]}}) }(A)$, then $(x,p(1))\in Z^1_{\rm sc}(\exp \g^{\Delta_{[0,1]}})$.
We have to show that if two elements $\eta_0=(x_0,e^{p_0(t)}*\partial_{0,1}x_0)$ and
 $\eta_1=(x_1,e^{p_1(t)}*\partial_{0,1}x_1)$ in $\MC_{{\Tot}_{TW}(\g^{\Delta_{[0,1]}})}(A)$
are homotopy equivalent, then $\Phi_{[0,1]}(\eta_0)\sim \Phi_{[0,1]}(\eta_1)$ in
$Z^1_{\rm sc}(\exp \g^{\Delta_{[0,1]}})$. Let $z(\xi,d\xi)$ be an homotopy
between $\eta_0$ and $\eta_1$. Therefore, $z(\xi,d\xi)$ is a Maurer-Cartan element for ${\Tot}_{TW}(\g^{\Delta_{[0,1]}})[\xi,d\xi]$ and so, reasoning as in the proof of Proposition \ref{prop.mc[01]}, we find
\[
z(\xi,d\xi)= (e^{T(\xi)}* u, e^{U(t,dt;\xi)}* v),
\]
with $T(0)=U(t,dt;0)=0$.
Since   $z(0)=\eta_0$, we get
\[
z(\xi,d\xi)=(e^{T(\xi)}* x_0, e^{U(t,dt; \xi)}* e^{p_0(t)}*
 \partial_{0,1}x_0).
\]
The face conditions for $z(\xi,d\xi)$ and uniqueness imply
\[ U(0;\xi)=\partial_{0,1}T(\xi) \quad \mbox{ \and }
\quad U(1;\xi)=\partial_{1,1}T(\xi). \] Moreover, $z(1)=\eta_1$, and so
\[ (e^{T(1)}* x_0, e^{U(t,dt; 1)}* e^{p_0(t)}* \partial_{0,1}x_0)=(x_1,e^{p_1(t)}*\partial_{0,1}x_1);   \]
by uniqueness again, we have
\[ e^{T(1)}*x_0=x_1.\]
Furthermore, $e^{U(t,dt; 1)}* e^{p_0(t)}* \partial_{0,1}x_0 = e^{p_1(t)}*\partial_{0,1}x_1$,
so, using the face conditions for $\eta_0$ and $\eta_1$, we obtain
\[
\partial_{0,1} x_0 = e^{-p_0(t)
 \bullet
 -U(t,dt; 1) \bullet p_1(t) \bullet \partial_{0,1} T(1)} *
\partial_{0,1} x_0\]
Next, we recall \cite[Lemma 6.15]{Dona} that if $L$ is a DGLA, $x(t,dt)$ is a Maurer-Cartan element for $L[t,dt]$ and $\mu(t,dt)\in L[t,dt]^0$ is such that $e^{\mu(t,dt)}*x(t,dt)=x(t,dt)$, then $\mu(1)$ is an element of the irrelevant stabilizer of $x(1)$. Therefore, in our case we get
\[-p_0(1)\bullet -\partial_{1,1}T(1) \bullet p_1(1) \bullet
\partial_{0,1} T(1) \in {\rm Stab}(\partial_{0,1} x_0).\]
\end{proof}

\begin{proposition} \label{prop.iso caso [0,1]}
Let $\g^{\Delta}$ be a semicosimplicial DGLA. The map
\[
\Phi_{[0,1]}:\Def_{ \Tot_{TW}(\g^{\Delta_{[0,1]}}) }  \to
H^1_{\rm sc}(\exp \g^{\Delta_{[0,1]}}) \]
 is an isomorphism of  functors of Artin rings. In particular, $H^1_{\rm sc}(\exp \g^{\Delta_{[0,1]}})$
 is a deformation functor.
\end{proposition}
\begin{proof}
Let $\Psi_{[0,1]}: Z^1_{\rm sc}(\exp
\g^{\Delta_{[0,1]}})(A) \to  \Tot_{TW}(\g^{\Delta_{[0,1]}})
\otimes{\mathfrak m}_A $ be the map given by $ (l, m) \mapsto (l, e^{t m}* \partial_{0,1}l)$; it is immediate to check that $\Phi_{[0,1]}$ actually takes its values in $\MC_{ \Tot_{TW}(\g^{\Delta_{[0,1]}}) }(A)$. Moreover, $\Psi_{[0,1]}$ induces a map \[
H^1_{\rm sc}(\exp \g^{\Delta_{[0,1]}})(A)\to \Def_{ \Tot_{TW}(\g^{\Delta_{[0,1]}}) }(A),
\]
which is the inverse of $\Phi_{[0,1]}$.
 Indeed, if $(l_0,m_0)\sim(l_1,m_1)$ in $Z^1_{\rm sc}(\exp
\g^{\Delta_{[0,1]}})(A)$, then
 there exist elements $a \in \g^0_0\otimes \m_A$ and $b\in{\mathfrak g}_1^{-1}\otimes{\mathfrak m}_A$ such that
\[
\begin{cases}
e^a * l_0=l_1\\
- m_0\bullet -\partial_{1,1}a \bullet m_1
\bullet \partial_{0,1}a=db+[\partial_{0,1}l_0,b].
\end{cases}
\]
Therefore, the images $(l_0, e^{t m_0}* \partial_{0,1}l_0)$ and  $(l_1, e^{t
m_1}* \partial_{0,1}l_1)$ are homotopic via the element
\[ z(\xi, d\xi)= ( e^{\xi a}* l_0, e^{t\bigl(\partial_{1,1}(\xi a)
\bullet m_0 \bullet (d(\xi b) +
[\partial_{0,1}l_0, \xi b]) \bullet -\partial_{0,1}(\xi a)\bigr)\bullet
\partial_{0,1}(\xi a)}* \partial_{0,1} l_0 ).\]
The composition $\Phi_{[0,1]}\circ \Psi_{[0,1]}\colon Z^1_{\rm sc}(\exp
\g^{\Delta_{[0,1]}})(A)\to Z^1_{\rm sc}(\exp
\g^{\Delta_{[0,1]}})(A)$ is clearly the
identity, whereas the composition
$\Psi_{[0,1]} \circ \Phi_{[0,1]}:  \MC_{
\Tot_{TW}(\g^{\Delta_{[0,1]}}) }(A) \to \MC_{ \Tot_{TW}(\g^{\Delta_{[0,1]}})
}(A)$ is homotopic to the identity. Indeed, $(x, e^{p(t)}* \partial_{0,1}x)$ and
$(x, e^{t p(1)}* \partial_{0,1}x)$ are homotopic in $\MC_{
\Tot_{TW}(\g^{\Delta_{[0,1]}}) }(A)$ via the element
$z(\xi,d\xi)= (x, e^{\xi t p(1) + (1-\xi) p(t) }*
\partial_{0,1}x)$.
\end{proof}

\begin{remark}A particular case of Proposition \ref{prop.iso caso [0,1]}, with an almost identical proof, has been considered by one of the authors in \cite{Dona}. Namely, given three DGLAs $L,M$ and $N$ and two DGLA morphisms $h\colon L\to M$ and $g\colon N\to M$, one can consider the semicosimplicial DGLA
  \[
\xymatrix{ L\oplus N
\ar@<2pt>[r]^{(0,g)}\ar@<-2pt>[r]_{(h,0)} & { M}
      \ar@<4pt>[r] \ar[r] \ar@<-4pt>[r] & 0
\ar@<6pt>[r] \ar@<2pt>[r] \ar@<-2pt>[r] \ar@<-6pt>[r]&
\cdots}
\]
to reobtain \cite[Theorem 6.17]{Dona}.
\end{remark}

\section{Proof of the main theorem}\label{sec.main}
In this section, we prove the existence of a  natural isomorphism
of functors of Artin rings $H^1_{\rm sc}(\exp \g^\Delta) \cong
\Def_{ \Tot_{TW}(\g^{ \Delta_{[0,2]} })}$, for  any
semicosimplicial DGLA $\g^{\Delta}$ such that $H^{-1}(\g_2)=0$.
As an immediate consequence we obtain a natural isomorphism of deformation
functors $H^1_{\rm sc}(\exp \g^\Delta) \cong \Def_{
{\Tot_{TW}}(\g^{ \Delta})}$, for any semicosimplicial DGLA
$\g^{\Delta}$, such that $H^{j}(\g_i)=0$ for $i\geq 0$ and $j<0$.
\par
The proof is considerably harder than in the case $\g^{\Delta_{[0,1]}}$ considered in the previous section. Indeed, we are still able to define a map $\Phi\colon \MC_{\Tot_{TW}(\mathfrak{g}^{\Delta_{[0,2]}})}
\to Z^1_{\rm sc}(\exp \g^\Delta)$ inducing a natural transformation  $\Def_{\Tot_{TW}(\mathfrak{g}^{\Delta_{[0,2]}})}
\to H^1_{\rm sc}(\exp \g^\Delta)$, but we will not be able to explicitly define an homotopy inverse to $\Phi$, so we will have to directly check that the map $\Def_{\Tot_{TW}(\mathfrak{g}^{\Delta_{[0,2]}})}
\to H^1_{\rm sc}(\exp \g^\Delta)$ is an isomorphism.

\begin{proposition} \label{prop.mappa phi}
Let $\g^\Delta$ be a semicosimplicial DGLA. The map
\[
\Phi:\MC_{ \Tot_{TW}(\g^{\Delta_{[0,2]}}) }(A)
\to (\g_0^1\oplus \g_1^0  ) \otimes \mathfrak m_A,
\]
given by
\[
(x, e^{p(t)}* \partial_{0,1}x, e^{q(s_0,s_1)+ r(s_0,s_1,ds_1,ds_1)}*\partial_{0,2}\partial_{0,1}x)
\mapsto (x, p(1) ),
\]
induces a natural transformation of functors of Artin rings
\[ \Def_{\Tot_{TW}(\mathfrak{g}^{\Delta_{[0,2]}})}
\to H^1_{\rm sc}(\exp \g^\Delta).\]
\end{proposition}

\begin{proof}
First we check that $\Phi$ takes its values in $Z^1_{\rm sc}(\exp \g^\Delta)(A)$.
The only nontrivial point consists in showing that  $- \partial_{2,2}p(1) \bullet \partial_{1,2}p(1)
\bullet -\partial_{0,2}p(1)$ is an element of the irrelevant stabilizer of $\partial_{2,2}\partial_{0,1}x$. This follows by the face condition
\[
e^{(- \partial_{2,2}p(t))\bullet(q(t,1-t)+ r(t,1-t,dt))\bullet(-q(0,1))}*\partial_{2,2}\partial_{0,1}x=\partial_{2,2}\partial_{0,1}x,
\]
applying \cite[Lemma 6.15]{Dona} once again. Next, we notice that
the  equivalence relation $\sim$ on $Z^1_{\rm sc}(\exp
\g^\Delta)(A)$ only involves the DGLAs ${\mathfrak g}_0$ and
${\mathfrak g}_1$; hence, we can conclude verbatim following the
proof of Proposition \ref{prop.mappa phi1}.

\end{proof}

\begin{proposition} \label{prop.surjectivity}
The map $\Phi: \Def_{{\Tot}_{TW}(\g^{\Delta_{[0,2]}})}(A)
\to H^1_{\rm sc}(\exp \g^\Delta)(A)$ is surjective.
\end{proposition}
\begin{proof}
Let $(l,m)\in Z^1_{\rm sc}(\exp \g^\Delta)(A)$ and $n\in \g_2^{-1}\otimes \m_A$, such that $ \partial_{0,2}m
\bullet - \partial_{1,2}m \bullet  \partial_{2,2}m = dn +
\frac{1}{2} [\partial_{2,2}\partial_{0,1}l, n]$.
Consider the element $w(t) = d(tn) + \frac{1}{2} [\partial_{2,2}
\partial_{0,1}l,t n]$ in the irrelevant stabilizer of $\partial_{2,2}\partial_{0,1}l$ and
\[ R(s_0,s_1)=s_0s_1\frac{s_0\partial_{2,2}m  \bullet -w(s_0) \bullet s_0
\partial_{0,2}m \bullet -s_0 \partial_{1,2}m}{s_0 (1-s_0) }
  \bullet s_0 \partial_{1,2}m \bullet s_1 \partial_{0,2}m. \]
Then,
\[ (l, e^{tm}* \partial_{0,1}l,  e^{R(s_0,s_1)}*
\partial_{0,2}\partial_{0,1}l) \]
is an element in $\MC_{\Tot_{TW}(\g^{\Delta_{[0,2]}})}(A)$ in the
fiber of $\Phi$ over $(l,m)$. Indeed, clearly it satisfies the Maurer-Cartan equation in $\g_0 \oplus \g_1\otimes \Omega_1 \oplus \g_2\otimes \Omega_2$; the first face
conditions follow easly noticing that
$R(0 ,t )= t \partial_{0,2}m$ and $R(t,0)= t
 \partial_{1,2}m$;
for the last one, we have:
\[ e^{R(t,1-t)}*\partial_{0,2}\partial_{0,1}l =
e^{ t\partial_{2,2}m \bullet -w(t) \bullet
\partial_{0,2}m} *\partial_{0,2}\partial_{0,1}l=
 \]
 \[
 = e^{ t\partial_{2,2}m    \bullet -w(t)} *
 \partial_{0,2}\partial_{1,1}l =
   e^{ t\partial_{2,2}m \bullet -w(t)} * \partial_{2,2}
   \partial_{0,1}l = e^{ t\partial_{2,2}m} *
 \partial_{2,2}\partial_{0,1}l. \]
\end{proof}

We will prove that the map
$\Phi:\Def_{\Tot_{TW}(\g^{\Delta_{[0,2]}})} (A) \to H^1_{\rm sc}(\exp
\g^\Delta)(A)$ is injective, under the hypothesis $H^{-1}(\g_2)=0$. For this we need two remarks.

\begin{remark} \label{rmk. defN modif banale}
Let $(L,d,[\ , \ ])$ be a DGLA, $A\in \bf{Art}_{\K}$ and
$x\in L^1\otimes \m_A$. The linear endomorphism $d_x=d+[x , \ ]$ of
$L\otimes \m_A$ is a differential if and only if $x\in\MC_L(A)$, and in this case
$(L\otimes \m_A, d_x, [\ , \ ])$ is a DGLA. So,
we can define the set of the Maurer-Cartan elements
$\MC_{L}^{x}(A)$ and the gauge action of $(L^0\otimes {\mathfrak m}_A, d_x, [\ , \ ])$ on it. We denote by $\Def_L^x(A)$ the quotient of
$\MC_L^x(A)$ with respect to the gauge action. The affine map
$$ \begin{array}{rll}
L \otimes \m_A & \to & L \otimes \m_A \\
v & \mapsto & v-x.
\end{array}$$
induces an isomorphism $\Def_L(A)\cong \Def_L^x(A)$ with obvious inverse $v\mapsto v+x$.

Next,  let $M\subseteq L$ be a sub-DGLA and let $x \in \MC_L(A)$. If
$M\otimes \m_A$ is closed under the differential $d_x$, then we can consider the set of Maurer-Cartan elements
$\MC_{M}^{x}(A)$, and its  quotient $\Def_M^x(A)$. The tangent space to $\Def_M^x(A)$ is $H^1(M\otimes {\mathfrak m}_A,d_x)$; so, by upper semicontinuity of cohomology, $H^1(M,d)=0$ implies that $\Def^x_{M}(A)$ is trivial, for all $x\in \MC_L(A)$ such that $d_x(M\otimes \m_A)\subseteq
M\otimes \m_A$.
\end{remark}

\begin{remark} \label{rmk.tr sur}
For any semicosimplicial DGLA $\g^{\Delta}$,  the truncation morphism
\[
\Tot^0_{TW}(\g^{\Delta_{[0,2]}}) \to
\Tot^0_{TW}(\g^{\Delta_{[0,1]}})
\] is surjective, i.e., for any  $(a_0,a_1) \in \Tot^0_{TW}(\g^{\Delta_{[0,1]}})$ there exist $a_2\in(\g_2\otimes\Omega_2)^0$ such that $(a_0,a_1,a_2)\in   \Tot^0_{TW}(\g^{\Delta_{[0,2]}})$. To see this, write $a_1(t,dt)= a_1^0(t)+ a_1^{-1}(t)dt$; then a possible choice for $a_2$ is
\[
a_2(s_0,s_1,ds_0,ds_1)= a^0_2(s_0,s_1)+ a^{-1}_{2,0}(s_0,s_1)ds_0 +a^{-1}_{2,1}(s_0,s_1)ds_1+ a^{-2}_{2}(s_0,s_1)ds_0ds_1,
\]
with
\begin{align*}a^0_2(s_0,s_1)&=\partial_{1,2}a^0_1(s_0) + \partial_{0,2}a^0_1(s_1)  -
\partial_{1,2}a^0_1(0)  \\
&\qquad\qquad+s_1 \frac{\partial_{2,2}a^0_{1}(s_0)-
\partial_{1,2}a^0_{1}(s_0)-\partial_{0,2}a^0_{1}(1-s_0)+\partial_{0,2}a^0_{1}(0)}{1-s_0};\\ &\\
 a^{-1}_{2,0}(s_0,s_1)&=\partial_{1,2}a_1^{-1}(s_0) +\frac{s_1}{1-s_0}
 \biggl(\partial_{2,2}a_1^{-1}(s_0) -
\partial_{1,2}a_1^{-1}(s_0) \\
&\qquad\qquad\qquad\qquad\qquad\qquad +
\partial_{0,2}a_1^{-1}(1-s_0) - s_0 \partial_{0,2}a_1^{-1}(0)\biggr );\\ & \\
a^{-1}_{2,1}(s_0,s_1)&= \partial_{0,2}a_1^{-1}(s_1) ds_1 -s_0\partial_{0,2}a_1^{-1}(0); \\ &\\
a^{-2}_{2}(s_0,s_1)&=0.
\end{align*}
It is an easy computation to verify that the element $(a_0,a_1,a_2)$ actually satisfies the face conditions.
\end{remark}

\begin{proposition} \label{prop.injectivity}
Let $\g^{\Delta}$ be a semicosimplicial DGLA, such that
$H^{-1}(\g_2)=0$. The map $\Phi:
\Def_{{\Tot}_{TW}(\g^{\Delta_{[0,2]}})}(A) \to H^1_{\rm sc}(\exp
\g^\Delta)(A)$ is injective.
\end{proposition}
\begin{proof}
Consider the commutative diagram
$$
\xymatrix{\Def_{\Tot_{TW}(\g^{\Delta_{[0,2]}})}(A)
\ar[r]^-{\Id_{[0,1]}} \ar[d]^-{\Phi}
&\Def_{\Tot_{TW}(\g^{\Delta_{[0,1]}})}(A)
\ar[d]_{\cong}^-{\Phi_{[0,1]}} \\ H^1_{\rm sc}(\exp \g^\Delta) (A)
\ar[r]^-{i} &H^1_{\rm sc}(\exp \g^{\Delta_{[0,1]}})(A); }
$$
since the map $\Phi_{[0,1]}$ is an isomorphism by Proposition
\ref{prop.iso caso [0,1]}, it is sufficient to prove that $\Id_{[0,1]}$ is injective.
Let $(x_0,x_1,x_2)$ and $(x'_0, x'_1, x'_2)$ be two Maurer-Cartan elements for $\Tot_{TW}(\g^{\Delta_{[0,2]}})$, such that $(x_0,x_1)$
and $(x'_0, x'_1)$ are gauge equivalent elements in
$\MC_{\Tot_{TW}(\g^{\Delta_{[0,1]}})}(A)$. Let $(a_0,a_1) \in \Tot^0_{TW}(\g^{\Delta_{[0,1]}})\otimes \m_A$ be an element realizing the gauge equivalence between $(x'_0,x'_1)$
and $(x_0, x_1)$, and let $(a_0,a_1,a_2)$ be a lift of $(a_0,a_1)$ in
$\Tot^0_{TW}(\g^{\Delta_{[0,2]}})\otimes \m_A$ (see Remark \ref{rmk.tr sur}). Then $(x'_0,x'_1,x'_2)$ is gauge equivalent via $(a_0,a_1,a_2)$ to the Maurer-Cartan element $(x_0,x_1,e^{a_2}*x'_2)$ and we are left to prove
that $(x_0,x_1,e^{a_2}*x'_2)$ is gauge equivalent to $(x_0,x_1,x_2)$.

To see this, consider the DGLA $\Tot_{TW}(\g^{\Delta_{[0,2]}})\otimes \m_A$ and modify its
differential with the Maurer-Cartan element $(x_0,x_1,x_2)$, as in  Remark \ref{rmk. defN modif banale}. Translation by $(x_0,x_1,x_2)$ gives an isomorphism
\[
\Def_{\Tot_{TW}
(\g^{\Delta_{[0,2]}})}(A)\cong \Def^{(x_0,x_1,x_2)}_{\Tot_{TW}(\g^{\Delta_{[0,2]}})}(A);
\]
hence $(x_0,x_1,x_2)$ and $(x_0,x_1,e^{a_2}*x'_2)$ will be gauge equivalent in $\MC_{\Tot_{TW}
(\g^{\Delta_{[0,2]}})}(A)$ if and only if $(0,0,0)
$ and $(0,0,e^{a_2}*x'_2-x_2)$ are gauge-equivalent in
$\MC^{(x_0,x_1,x_2)}_{\Tot_{TW}(\g^{\Delta_{[0,2]}})}(A)$.
Next, observe that the sub-DGLA
$\Tot_{TW}(\g^{\Delta_{[2,2]}})\otimes \m_A$ of  $\Tot_{TW}(\g^{\Delta_{[0,2]}})\otimes
\m_A$ is closed under the
modified differential $d_{(x_0,x_1,x_2)}$, so we can consider the deformation functor $\Def^{(x_0,x_1,x_2)}_{\Tot_{TW}
(\g^{\Delta_{[2,2]}})}(A)$. Since $H^{1}(\Tot_{TW}
(\g^{\Delta_{[2,2]}}),d_{TW})=H^1(\Tot
(\g^{\Delta_{[2,2]}}),d_{\rm Tot})=H^{-1}(\g_2)=0$, this deformation functor is trivial (see Remark \ref{rmk. defN modif banale}). Therefore $(0,0,e^{a_2}*x'_2-x_2)$ is gauge equivalent to $(0,0,0)$ as an element of $\MC^{(x_0,x_1,x_2)}_{\Tot_{TW}
(\g^{\Delta_{[2,2]}})}(A)$, and so, a fortiori, as an element of $\MC^{(x_0,x_1,x_2)}_{\Tot_{TW}
(\g^{\Delta_{[0,2]}})}(A)$.
\end{proof}
Summing up, and recalling Proposition \ref{prop.tre suff}, we have proved:
\begin{theorem} \label{th.iso funt}
Let $\g^{\Delta}$ be a semicosimplicial DGLA,  and let
$\Tot_{TW}(\g^{\Delta})$ and  $\Tot_{TW}(\g^{\Delta_{[0,2]}})$ be
the Thom-Whitney DGLAs associated with $\g^{\Delta}$ and
$\g^{\Delta_{[0,2]}}$, respectively. Assume that $H^{-1}(\g_2)=0$;
then, there is a natural isomorphism of funtors
$\Def_{\Tot_{TW}(\g^{\Delta_{[0,2]}})} \cong H^1_{\rm sc}(\exp
\g^\Delta)$. If moreover $H^{j}(\g_i)=0$ for all $i\geq 0$ and
$j<0$, then there is a natural isomorphism of funtors
$\Def_{\Tot_{TW}(\g^{\Delta})} \cong H^1_{\rm sc}(\exp
\g^\Delta)$. In particular, in this case, the tangent space to
$H^1_{\rm sc}(\exp \g^\Delta)$ is $H^1({\rm Tot}(\g^\Delta))$ and
obstructions are contained in $H^2({\rm Tot}(\g^\Delta))$.
\end{theorem}

\begin{theorem}
Let $X$ be a paracompact Hausdorff topological space, and let
$\mathcal L$ be a sheaf of differential graded Lie algebras on $X$, such that the DGLAs $\mathcal L(U_{i_0  \ldots  i_k})$ has no negative cohomology. Then, every
refinement ${\mathcal V}\geq {\mathcal U}$ of open covers of
$X$ induces a natural morphism of deformation functors
$\Def_{\Tot_{TW}({\mathcal L}({\mathcal U}))}\to
\Def_{\Tot_{TW}({\mathcal L}({\mathcal V}))}$. In
particular, the direct limit
\[
\Def_{[{\mathcal L}]} = \lim_{\stackrel{\longrightarrow}{\mathcal U}}
 \Def_{\Tot_{TW}({\mathcal L}({\mathcal
U}))}
\]
is well defined and there is 
natural isomorphism of functors of Artin rings
\[
H_{\rm Ho}^1(X;\exp{\mathcal L})\cong \Def_{[{\mathcal L}]}.
\]
Moreover, if acyclic open covers for $\mathcal L$ are cofinal
in the directed family of all open covers of $X$, then \[ H_{\rm
Ho}^1(X; \exp \mathcal L) \cong H_{\rm sc}^1(\exp\mathcal
L(\U))\qquad \text{and}\qquad \Def_{[{\mathcal L}]}\cong
\Def_{\Tot_{TW}({\mathcal L}({\mathcal U}))},
\]
for every ${\mathcal L}$-acyclic open cover ${\mathcal U}$ of $X$.
\end{theorem}

\begin{proof}
Let ${\mathcal V}\geq {\mathcal U}$ be a refinement of open
covers of $X$, and let $\tau$ be a refinement function, it induces a natural morphism of semicosimplicial Lie algebras
${\mathcal L}({\mathcal U})\to {\mathcal L}({\mathcal V})$ and so a commutative diagram of
natural transformations
\[\xymatrix{
\Def_{\Tot_{TW}({\mathcal L}({\mathcal
U}))}\ar[d]\ar[r]^-{\sim}& H^1_{\rm
sc}(\exp{\mathcal L}({\mathcal
U})) \ar[d]\\
\Def_{\Tot_{TW}({\mathcal L}({\mathcal
V}))}\ar[r]^-{\sim}&H^1_{\rm sc}(\exp{\mathcal
L}({\mathcal V})).}
\]
Horizontal arrows are isomorphisms by Theorem \ref{th.iso funt}, and the right vertical
arrow is independent of the refinement function $\tau$, as
observed in Lemma \ref{lemma limite raffinamenti}. Hence, also the
left morphism is independent of $\tau$, then the direct limit
\[
\Def_{[{\mathcal L}]} = \lim_{\stackrel{\longrightarrow}{\mathcal
U}} \Def_{\Tot_{TW}({\mathcal L}({\mathcal U}))}
\]
is well defined and  we have a  natural isomorphism
$\Def_{[{\mathcal L}]}\cong H_{\rm Ho}^1(X;\exp{\mathcal L})$.
Assume now that acyclic open covers for ${\mathcal L}$ are
cofinal in the family of all open covers of $X$. Then, for any
refinement ${\mathcal V}\geq {\mathcal U}$ of acyclic open
covers, the  DGLAs-morphism $\Tot_{TW}({\mathcal L}({\mathcal
U}))\to \Tot_{TW}({\mathcal L}({\mathcal V}))$ is a
quasi-isomorphism by Leray's theorem. Therefore, we have a
commutative diagram of natural transformations
\[\xymatrix{
\Def_{\Tot_{TW}({\mathcal L}({\mathcal
U}))}\ar[d]^\wr\ar[r]^-{\sim}&H^1_{\rm
sc}(\exp{\mathcal L}({\mathcal
U}))\ar[d]\\
\Def_{\Tot_{TW}({\mathcal L}({\mathcal
V}))}\ar[r]^-{\sim}&H^1_{\rm sc}(\exp{\mathcal
L}({\mathcal V})),\\ }
\]
where also the right vertical arrow is forced to be an
isomorphism. Taking the direct limit over ${\mathcal L}$-acyclic
covers, we obtain that, if ${\mathcal U}$ is an ${\mathcal
L}$-acyclic open cover of $X$, then $ H_{\rm Ho}^1(X; \exp
\mathcal L) \cong H^1_{\rm sc}(\exp\mathcal L(\mathcal U))$ and $
\Def_{[{\mathcal L}]}\cong \Def_{\Tot_{TW}({\mathcal L}({\mathcal
U}))}$.
\end{proof}

\section{Conclusions and further developements}
We can now sum up our results to obtain a DGLA description of
infinitesimal deformations of a coherent sheaf. In Section
\ref{sec.Example}, we analised infinitesimal deformations of a
coherent sheaf $\mathcal F$ of $\O_X$-modules on a ringed space
$(X,\O_X)$. If $\Eps^{\cdot} \to \mathcal F\to 0$ is
a locally free resolution of $\mathcal F$ on $X$, we showed how
infinitesimal deformations of ${\mathcal F}$ can be expressed in terms
of the sheaf of DGLAs $\Eps nd^* (\Eps^{\cdot})$. More precisely, in
Section \ref{sec.semicosimplicial dglas}, we showed that the functor
of infinitesimal deformations of ${\mathcal F}$ is isomorphic to
$H^1_{\rm Ho}(X; \exp {\Eps nd}^*(\Eps^{\cdot}))$.
\par
Since negative Ext-groups between coherent sheaves are always trivial,
all terms in the semicosimplicial DGLA $\Eps nd^*(\Eps^{\cdot})(\U)$
have zero negative cohomology. Therefore,
Theorem \ref{th.iso funt} applies and we obtain that the functor of
infinitesimal deformations of ${\mathcal F}$ is isomorphic to
$\Def_{[\Eps nd^*(\Eps^{\cdot})]}$; in particular, we recover the well
known
fact that the tangent space to $\Def_{\mathcal F}$ is
$\Ext^1({\mathcal F},{\mathcal F})$ and that its obstructions are
contained in $\Ext^2({\mathcal F},{\mathcal F})$.
\par
Moreover, if $X$ is a smooth complex variety, then the DGLA controlling
infinitesimal deformations of $\mathcal F$ turns out to be not at all
mysterious. Indeed, let $ {\Eps nd}^*(\Eps^{\cdot}) \to
\mathcal{A}^{0,*}_X(\Eps nd^*(\Eps^{\cdot})) $ be the Dolbeault
resolution of ${\Eps nd}^*(\Eps^{\cdot})$. Since this resolution is
fine, by Theorem \ref{th.semicosimpl vs dgla} the functor of
infinitesimal deformations of ${\mathcal F}$ is isomorphic to the
deformation functor associated with the DGLA $A^{0,*}_X(\Eps
nd^*(\Eps^{\cdot}))$ of global sections of $\mathcal{A}^{0,*}_X(\Eps
nd^*(\Eps^{\cdot}))$. We can also give an explicit description of this
isomorphism of deformation functors. Indeed, a natural isomorphism
\[ \Def_{A^{0,*}_X(\Eps nd^*(\Eps^{\cdot}))}(B) \to \Def_{\mathcal
F}(B), \qquad \mbox{for}\ B \in \bf{Art}_{\K}   \]
is defined by associating with every Maurer-Cartan element $\xi$ of
the DGLA $A^{0,*}_X(\Eps nd^*(\Eps^{\cdot}))$ the cohomology
sheaf of $({\mathcal A}_X^{0,*}(\Eps^{\cdot}) \otimes B, \deltabar
+d_{\Eps^{\cdot}}+ \xi)$. Note that, by semicontinuity, this cohomology sheaf is concentrated in degree zero.

\bigskip
The techniques developed in this paper apply to a wide range of other
geometric examples. More explicitly, we can use them in all cases
when local deformations admit
a simple DGLA description in terms of a resolution of the object to be
deformed, for instance, in the case of infinitesimal deformations of a
singular variety. Namely, let $X$ be a singular variety,  ${\mathcal
O}_X$  the sheaf of regular function of $X$ and  ${\mathcal
R}^\cdot\to {\mathcal O}_X$ its standard free resolution \cite[Section
1.5]{Illusie}. Then, the deformation functor of infinitesimal
deformations of $X$ is isomorphic to $ H_{\rm Ho}^1(X; \exp {{\mathcal
D}er}^*({\mathcal R}^\cdot) )$; see \cite{dde2} for details. From
this, we also recover the classical result that the tangent space to
deformations of $X$ is $\Ext^1({\mathbb L}_X,{\mathcal O}_X)$, and
that obstructions are contained in  $\Ext^2({\mathbb L}_X,{\mathcal
O}_X)$, where ${\mathbb L}_X$ is the cotangent complex of $X$.


\begin{thebibliography}{99}
\bibitem{Artin}
M.~Artin: \emph{Deformation of singularities.} Tata institute of foundamental research, Bombay, (1976).

\bibitem{chenggetzler} X.~Z.~Cheng, E.~Getzler:
\emph{Homotopy commutative algebraic structures.} J. Pure Appl.
Algebra, \textbf{212}, (2008), 2535-2542.
\bibitem{dupont1}
J.~L.~Dupont: \emph{Simplicial de Rham cohomology and
characteristic classes of flat bundles.}  Topology,  \textbf{15},
(1976), 233-245.

\bibitem{dupont2}
J.~L.~Dupont: \emph{Curvature and characteristic classes.}
Lecture Notes in Mathematics, \textbf{640}, Springer-Verlag, (1978).


\bibitem{dde2}
D.~Fiorenza, D.~Iacono, E.~Martinengo:
\emph{Infinitesimal deformations of singular varieties.} (in prepapartion).

\bibitem{Fio-Man}
D.~Fiorenza, M.~Manetti: \emph{$L_{\infty}$-structures on mapping
cones.} Algebra \& Number Theory, \textbf{1}, (2007), 301-330.


\bibitem{Fio-Man-Mart}
D.~Fiorenza, M.~Manetti, E.~Martinengo: \emph{Semicosimplicial
DGLAs in deformation theory}. 
\texttt{arxiv:math.AG/08030399}.

\bibitem{getzler}
E.~Getzler: \emph{Lie theory for nilpotent $L_{\infty}$-algebras.}
Ann. of Math.,  \textbf{170}, (1), (2009), 271-301.



\bibitem{hinich}
V.~Hinich: \emph{Descent of Deligne groupoids.}
Int. Math. Res. Notices,  (1997),  \textbf{5}, 223-239.

\bibitem{hirschowitz-simpson} 
A.~Hirschowitz, C.~Simpson: \emph{Descent pour les n-champs.} 
\texttt{arXiv:9807049v3}. 

\bibitem{Dona}
D.~Iacono:
\emph{$L_{\infty}$-algebras and deformations of
holomorphic maps.} Int. Math. Res. Notices, (2008), \textbf{8}, Art. ID rnn013, 36 pp.

\bibitem{Illusie} L.~Illusie:
\emph{Complexe cotangent et deformations I, II.} Lecture Notes in
Mathematics, {\bf  239, 283}, Springer-Verlag, New York/Berlin,
(1971-1972).

\bibitem{bib kodaira libro} K.~Kodaira:
\emph{Complex Manifolds and Deformation of Complex Structures.} Die
Grundlehren der mathematischen Wissenschaften in
Einzeldarstellungen, {\bf 283}, Springer-Verlag, New York/Berlin,
(1986).

\bibitem{bib Kodaira SpencerII} K.~Kodaira, D.~C.~Spencer:
\emph{On Deformations of Complex Analytic Structures, II.}  Ann. of
Math., \textbf{67} (2), (1958), 403-466.




\bibitem{bib kuranishi} M.~Kuranishi:
\emph{Deformations of compact complex manifolds.} S\'eminaire de
Mathematiques Sup\'erieures, No. \textbf{39}, (\'Et\'e 1969), Les
Presses de l'Universit\'e de Montreal, Montreal, (1971).

\bibitem{lurie} J.~Lurie: \emph{Higher topos theory.} Annals of Mathematics Studies, \textbf{170}, Princeton University Press, Princeton, NJ, (2009).


\bibitem{Man Pisa} M.~Manetti:
\emph{Deformation theory via differential graded Lie algebras.}
Seminari di Geometria Algebrica 1998-1999,
Scuola Normale Superiore (1999).



\bibitem{bib Manetti IMNR} M.~Manetti:
\emph{Extended deformation functors.} Int. Math. Res. Not., 
{\bf14}, (2002), 719-756.


\bibitem{navarro}
V.~Navarro Aznar: \emph{Sur la th\'eorie de Hodge-Deligne.}
Invent. Math., \textbf{90}, (1987),11-76.



\bibitem{Pridham} J.~P.~Pridham:
\emph{Deformations via Simplicial Deformation Complexes.} 
\texttt{arXiv:math/0311168v6}.


\bibitem{quillen} D.~Quillen:
 \emph{Rational homotopy theory.}
 Ann. of Math., (2), \textbf{90}, (1969), 205-295.

\bibitem{Seidel.Thomas} P.~Seidel, R.~P.~Thomas:
\emph{Braid group actions on derived categories of coherent
sheaves.} Duke Math. J. \textbf{108}         
   (2001), no. 1, 37-108.

\bibitem{Sernesi} E.~Sernesi: \emph{Deformation of algebraic schemes.} Springer, \textbf{334}, (2006).

\bibitem{toen} B.~Toen: \emph{Higher and derived stack: a global overview.}
Proc. Sympos. Pure Math., \textbf{80}, Part 1, Amer. Math. Soc., Providence, RI, (2009).

\bibitem{yakutieli}
A.~Yekutieli: \emph{Twisted deformation quatization of algebraic varieties.}  \texttt{arXiv:0905.0488v2 }

\end{thebibliography}
\end{document}